\documentclass[11pt,a4paper]{article}

\usepackage{amsmath,amssymb}
\newcommand{\e}{\varepsilon}
\newcommand{\va}{\varphi}
\newcommand{\D}{\Delta}

\newcommand{\n}{\nabla}
\newcommand{\N}{\frac{N}{2}}
\newcommand{\NN}{\frac{N}{p}}

\newcommand{\p}{\partial}

\newcommand{\R}{\mathbb{R}}
\newcommand{\h}{\hookrightarrow}

\newtheorem{definition}{Definition}
\newtheorem{theorem}{Theorem}
\newtheorem{proposition}{Proposition}
\newtheorem{corollaire}{Corollary}

\newtheorem{notation}{Notation}
\newtheorem{remarka}{Remark}

\newtheorem{lemme}{Lemma}

\title{Well-posedness in critical spaces for barotropic viscous fluids}
\author{Boris Haspot \thanks{Karls Ruprecht Universit\"at Heidelberg, Institut for Applied Mathematics,
Im Neuenheimer Feld 294,
D-69120 Heildelberg, Germany.
Tel. 49(0)6221-54-6112 }}
\date{}
\begin{document}
\maketitle
\begin{abstract}
This paper is dedicated to the study of viscous compressible barotropic fluids in dimension $N\geq2$. We address
the question of well-posedness for {\it large} data having critical Besov regularity. 
Our result improve the analysis of R. Danchin in \cite{DW}, by the fact that we choose initial density more general in $B^{\NN}_{p,1}$ with $1\leq p<+\infty$.
Our result relies on a new a priori estimate for the velocity, where we introduce a new structure to \textit{kill} the coupling between the density and the velocity. In particular our result is the first where we obtain uniqueness without imposing hypothesis on the gradient of the density.
\end{abstract}
\section{Introduction}
The motion of a general barotropic compressible fluid is described by the following system:
\begin{equation}
\begin{cases}
\begin{aligned}
&\p_{t}\rho+{\rm div}(\rho u)=0,\\
&\p_{t}(\rho u)+{\rm div}(\rho u\otimes u)-{\rm div}(\mu(\rho)D(u))-\n(\lambda(\rho){\rm div} u)+\n P(\rho)=\rho f,\\
&(\rho,u)_{/t=0}=(\rho_{0},u_{0}).
\end{aligned}
\end{cases}
\label{0.1}
\end{equation}
Here $u=u(t,x)\in\R^{N}$ stands for the velocity field and $\rho=\rho(t,x)\in\R^{+}$ is the density.
The pressure $P$ is a suitable smooth function of $\rho$.
We denote by $\lambda$ and $\mu$ the two viscosity coefficients of the fluid,
which are assumed to satisfy $\mu>0$ and $\lambda+2\mu>0$ (in the sequel to simplify the calculus we will assume the viscosity coefficients as constants). Such a conditions ensures
ellipticity for the momentum equation and is satisfied in the physical cases where $\lambda+\frac{2\mu}{N}>0$.
We supplement the problem with initial condition $(\rho_{0},u_{0})$ and an outer force $f$.
Throughout the paper, we assume that the space variable $x\in\R^{N}$ or to the periodic
box ${\cal T}^{N}_{a}$ with period $a_{i}$, in the i-th direction. We restrict ourselves the case $N\geq2$.\\
The problem of existence of global solution in
time for Navier-Stokes equations was addressed in one dimension for
smooth enough data by Kazhikov and Shelukin in \cite{5K}, and for
discontinuous ones, but still with densities away from zero, by Serre
in  \cite{5S} and Hoff in \cite{5H1}. Those results have been
generalized to higher dimension by Matsumura and Nishida in
\cite{MN} for smooth data close
to equilibrium and by Hoff in the case of discontinuous data in \cite{5H2,5H3}. All
those results do not require to be far from the vacuum.
The existence and uniqueness of local classical solutions for (\ref{0.1})
with smooth initial data such that the density $\rho_{0}$ is bounded
and bounded away from zero (i.e.,
$0<\underline{\rho}\leq\rho_{0}\leq M$)
has been stated by Nash in \cite{Nash}. Let us emphasize that no stability condition was required there.
On the other hand, for small smooth perturbations of a stable
equilibrium with constant positive density, global well-posedness
has been proved in \cite{MN}. Many works on the case of the one dimension have been devoted
to the qualitative behavior of solutions for large time (see for
example \cite{5H1,5K}). Refined functional analysis has been used
for the last decades, ranging from Sobolev, Besov, Lorentz and
Triebel spaces to describe the regularity and long time behavior of
solutions to the compressible model \cite{5So}, \cite{5V},
\cite{5H4}, \cite{5K1}.
Let us recall that (local) existence and uniqueness for (\ref{0.1}) in the case of smooth
data with no vacuum has been stated for long in the pioneering works by J. Nash \cite{Nash},
and A. Matsumura, T. Nishida \cite{MN}. For results of weak-strong uniqueness, we refer to the work of P. Germain \cite{PG}.\\
Guided in our approach by numerous works dedicated to the incompressible Navier-Stokes equation (see e.g \cite{Meyer}):
$$
\begin{cases}
\begin{aligned}
&\p_{t}v+v\cdot\n v-\mu\D v+\n\Pi=0,\\
&{\rm div}v=0,
\end{aligned}
\end{cases}
\leqno{(NS)}
$$
we aim at solving (\ref{0.1}) in the case where the data $(\rho_{0},u_{0},f)$ have \textit{critical} regularity.\\
By critical, we mean that we want to solve the system functional spaces with norm
in invariant by the changes of scales which leaves (\ref{0.1}) invariant.
In the case of barotropic fluids, it is easy to see that the transformations:
\begin{equation}
(\rho(t,x),u(t,x))\longrightarrow (\rho(l^{2}t,lx),lu(l^{2}t,lx)),\;\;\;l\in\R,
\label{1}
\end{equation}
have that property, provided that the pressure term has been changed accordingly.\\
The use of critical functional frameworks led to several new weel-posedness results for compressible
fluids (see \cite{DL,DG,DW}). In addition to have a norm invariant by (\ref{1}),
appropriate functional space for solving (\ref{0.1}) must provide a control on the $L^{\infty}$
norm of the density (in order to avoid vacuum and loss of ellipticity). For that reason,
we restricted our study to the case where the initial data $(\rho_{0},u_{0})$ and external force $f$
are such that, for some positive constant $\bar{\rho}$:
$$(\rho_{0}-\bar{\rho})\in B^{\NN}_{p,1},\;u_{0}\in B^{\frac{N}{p_{1}}-1}_{p_{1},1}\;\;\mbox{and}\;\;f\in L^{1}_{loc}(\R^{+},\in B^{\frac{N}{p_{1}}-1}_{p_{1},1})$$
with $(p,p_{1})\in [1,+\infty[$ good chosen.
In \cite{DW}, however, we hand to have $p=p_{1}$, indeed in this article there exists a very strong coupling between the pressure
and the velocity. To be more precise, the pressure term is considered as a term of rest for the elliptic operator in the momentum equation of (\ref{0.1}). This paper improve the results of R. Danchin in \cite{DL,DW}, in the sense that the initial density belongs to larger spaces $B^{\NN}_{p,1}$ with $p\in[1,+\infty[$.
The main idea of this paper is to introduce a new variable than the velocity in the goal to \textit{kill}
the relation of coupling between the velocity and the density.
In the present paper, we address the question of local well-posedness in the critical functional framework under the assumption
that the initial density belongs to critical Besov space with a index of integrability different of this of the velocity.
We adapt the spirit of the results of \cite{AP} and \cite{H} which treat the case of Navier-Stokes incompressible with dependent density (at the difference than in these works the velocity and the density are naturally decoupled).
To simplify the notation, we assume from now on that $\bar{\rho}=1$. Hence as long as $\rho$ does not vanish, the equations for ($a=\rho^{-1}-1$,$u$) read:
\begin{equation}
\begin{cases}
\begin{aligned}
&\p_{t}a+u\cdot\n a=(1+a){\rm div}u,\\
&\p_{t}u+u\cdot\n u-(1+a){\cal A}u+\n (g(a))=f,
\end{aligned}
\end{cases}
\label{0.6}
\end{equation}
In the sequel we will note ${\cal A}=\mu\D+(\lambda+\mu)\n{\rm div}$ and where $g$ is a smooth function which may be computed from the pressure function $P$.\\
One can now state our main result.
\begin{theorem}
Let $P$ a suitably smooth function of the density and $1\leq p_{1}\leq p\leq 2N$ such that $\frac{1}{p_{1}}\leq\frac{1}{N}+\frac{1}{p}$
and $2\NN-1>0$. Assume that $u_{0}\in B^{\frac{N}{p_{1}}-1}_{p_{1},1}$, $f\in L^{1}_{loc}(\R^{+},B^{\frac{N}{p_{1}}-1}_{p_{1},1})$ and $a_{0}\in B^{\NN}_{p,1}$
 with
$1+a_{0}$ bounded away from zero. \\
If $\frac{1}{p}+\frac{1}{p_{1}}>\frac{1}{N}$ there exists a positive time $T$ such that system (\ref{0.1}) has a solution
$(a,u)$ with $1+a$ bounded away from zero,
$$a\in \widetilde{C}([0,T],B^{\NN}_{p,1}
),\;\;u\in \widetilde{C}([0,T];B^{\frac{N}{p_{1}}-1}_{p_{1},1}+B^{\NN+1}_{p,1})\cap\widetilde{L}^{1}(
B^{\frac{N}{p_{1}}+1}_{p_{1},1}+B^{\NN+2}_{p,1}).$$
Moreover this solution is unique if $\frac{2}{N}\leq\frac{1}{p}+\frac{1}{p_{1}}$.
\label{theo1}
\end{theorem}
\begin{remarka}
It seems possible to improve the theorem \ref{theo1} by choosing initial data $a_{0}$ in $B^{\NN}_{p,\infty}\cap B^{0}_{\infty,1}$, however some supplementary conditions appear on $p_{1}$ in this case.
\end{remarka}
The key to theorem \ref{theo1} is to introduce a new variable $v_{1}$ to control the velocity where to avoid the coupling between the density and the velocity, we analyze by a new way the pressure term. More precisely we write the gradient of the pressure as a Laplacian of the variable $v_{1}$, and we introduce this term in the linear part of the momentum equation. We have then a control on $v_{1}$ which can write roughly as $u-{\cal G}P(\rho)$ where ${\cal G}$ is a pseudodifferential operator of order $-1$. By this way, we have canceled the coupling between $v_{1}$ and the density, we next verify easily that we have a control Lipschitz of the gradient of $u$ (it is crucial to estimate the density by the transport equation).
\begin{remarka}
In the present paper we did not strive for unnecessary generality which may hide the new ideas of our analysis. Hence we focused on the somewhat academic model of barotropic fluids. In physical contexts however, a coupling with the energy equation has to be introduced. Besides, the viscosity coefficients may depend on the density. We believe that our analysis may be carried out to these more general models. (See \cite{H1}).
\end{remarka}
In \cite{5H5}, D. Hoff show a very strong theorem of uniqueness for the weak solution when the pressure is of the specific form
$P(\rho)=K\rho$ with $K>0$. Similarly in \cite{5H2}, \cite {5H3}, \cite{5H4}, D. Hoff get global weak solution with regularizing effects on the velocity. In particular when the pressure is on this form, he doe not need to have estimate on the gradient of the density.
In the following corollary, we will observe that this type of pressure assure a specific structure and avoid to impose that $p<2N$.
\begin{corollaire}
Assume that $P(\rho)=K\rho$ with $K>0$. Let $1\leq p_{1}\leq p\leq+\infty$ such that $\frac{1}{p_{1}}\leq\frac{1}{N}+\frac{1}{p}$. Assume that $u_{0}\in B^{\frac{N}{p_{1}}-1}_{p_{1},1}$, $f\in L^{1}_{loc}(\R^{+},B^{\frac{N}{p_{1}}-1}_{p_{1},1})$ and $a_{0}\in B^{\NN}_{p,1}$
with
$1+a_{0}$ bounded away from zero. 
If $\frac{1}{p}+\frac{1}{p_{1}}>\frac{1}{N}$ there exists a positive time $T$ such that system (\ref{0.1}) has a solution
$(a,u)$ with $1+a$ bounded away from zero,
$$a\in \widetilde{C}([0,T],B^{\NN}_{p,1})
,\;\;u\in \widetilde{C}([0,T];B^{\frac{N}{p_{1}}-1}_{p_{1},1}+B^{\NN+1}_{p,1})\cap\in \widetilde{L}^{1}(
B^{\frac{N}{p_{1}}+1}_{p_{1},1}+B^{\NN+1}_{p,1}).$$
If moreover we assume that $\sqrt{\rho_{0}}u_{0}\in L^{2}$, $\rho_{0}-\bar{\rho}\in L^{1}_{2}$, $u_{0}\in L^{\infty}$
and $\lambda=0$
then the solution $(a,u)$ is unique.
\label{coro11}
\end{corollaire}
\begin{remarka}
Here $L^{1}_{2}$ defines the corresponding Orlicz space.
\end{remarka}
\begin{remarka}
Up to my knowlledge, it seems that it is the first time that we get strong solution without condition of controll in space with positive regularity for the gradient of the density.
\end{remarka}
\begin{remarka}
Moreover we can observe that with this type of pressure we are very close to have existence of strong in finite time solution for initial data $(a_{0},u_{0})$ in
$B^{0}_{\infty,1}\times B^{1}_{N,1}$. It means that this theorem rely the result of D. Hoff where the initial density is assumed $L^{\infty}$ but where the initial velocity is more regular and the results of R. Danchin in \cite{DW}.
\end{remarka}
\begin{remarka}
In particular we can show that the solution of D. Hoff in \cite{5H4} are unique on a finite time interval $[0,T]$.
\end{remarka}
The study of the linearization of (\ref{0.1}) leads also the following continuation criterion:
\begin{theorem}
\label{theo3}
Let $1\leq p_{1}\leq p\leq+\infty$ such that $\frac{N}{p_{1}}-1\leq\frac{N}{p}$ and $\frac{N}{p_{1}}-1+\NN>0$. Assume that (\ref{0.1}) has a solution $(a,u)\in C([0,T),B^{\NN}_{p,1}\times(B^{\frac{N}{p_{1}}-1}_{p_{1},1}+B^{\NN+1}_{p,1})^{N})$ with $p_{1}>N$,
$\rho_{0}^{\frac{1}{p_{1}}}u_{0}\in L^{p_{1}}$ and:
\begin{equation}
\lambda\leq\frac{4\mu}{N^{2}(p_{1}-1)},
\label{inegaliteviscosite} 
\end{equation}
on the time
interval $[0,T)$ which satisfies the following three conditions:
\begin{itemize}
\item the function $a$ belongs to $L^{\infty}(0,T;B^{\NN}_{p,1}),$
\item the function $1+a$ is bounded away from zero.
\end{itemize}
Then $(a,u)$ may be continued beyond $T$.
\end{theorem}
\begin{remarka}
 Up my knowledge, it is the first time that we get a criterion of blow-up for strong solution for compressible Navier-Stokes system without imposing a controll Lipschitz of the norm $\n u$.
\end{remarka}
Our paper is structured as follows. In section \label{section2}, we give a few notation and briefly introduce the basic Fourier analysis
techniques needed to prove our result. In section
section \ref{section3} and \ref{section4} are devoted to the proof of key estimates for the linearized
system (\ref{0.1}). In section \ref{section5}, we prove the theorem \ref{theo1} and corollary \ref{theo11} whereas section \ref{section6}
is devoted to the proof of continuation criterions of theorem \ref{theo2} and \ref{theo3}. Two inescapable technical commutator estimates and some theorems of ellipticity are postponed in an appendix.
\section{Littlewood-Paley theory and Besov spaces}
\label{section2}
Throughout the paper, $C$ stands for a constant whose exact meaning depends on the context. The notation $A\lesssim B$ means
that $A\leq CB$.
For all Banach space $X$, we denote by $C([0,T],X)$ the set of continuous functions on $[0,T]$ with values in $X$.
For $p\in[1,+\infty]$, the notation $L^{p}(0,T,X)$ or $L^{p}_{T}(X)$ stands for the set of measurable functions on $(0,T)$
with values in $X$ such that $t\rightarrow\|f(t)\|_{X}$ belongs to $L^{p}(0,T)$.
Littlewood-Paley decomposition  corresponds to a dyadic
decomposition  of the space in Fourier variables.
Let $\alpha>1$ and $(\va,\chi)$ be a couple of smooth functions valued in $[0,1]$, such that $\va$ is supported in the shell
supported in
$\{\xi\in\R^{N}/\alpha^{-1}\leq|\xi|\leq2\alpha\}$, $\chi$ is supported in the ball $\{\xi\in\R^{N}/|\xi|\leq\alpha\}$
such that:
$$\forall\xi\in\R^{N},\;\;\;\chi(\xi)+\sum_{l\in\mathbb{N}}\varphi(2^{-l}\xi)=1.$$
Denoting $h={\cal{F}}^{-1}\varphi$, we then define the dyadic
blocks by:
$$
\begin{aligned}
&\D_{l}u=0\;\;\;\mbox{if}\;\;l\leq-2,\\
&\D_{-1}u=\chi(D)u=\widetilde{h}*u\;\;\;\mbox{with}\;\;\widetilde{h}={\cal F}^{-1}\chi,\\
&\D_{l}u=\varphi(2^{-l}D)u=2^{lN}\int_{\R^{N}}h(2^{l}y)u(x-y)dy\;\;\;\mbox{with}\;\;h={\cal F}^{-1}\chi,\;\;\mbox{if}\;\;l\geq0,\\
&S_{l}u=\sum_{k\leq
l-1}\D_{k}u\,.
\end{aligned}
$$
Formally, one can write that:
$u=\sum_{k\in\mathbb{Z}}\D_{k}u\,.$
This decomposition is called nonhomogeneous Littlewood-Paley
decomposition. 
\subsection{Nonhomogeneous Besov spaces and first properties}
\begin{definition}
For
$s\in\R,\,\,p\in[1,+\infty],\,\,q\in[1,+\infty],\,\,\mbox{and}\,\,u\in{\cal{S}}^{'}(\R^{N})$
we set:
$$\|u\|_{B^{s}_{p,q}}=(\sum_{l\in\mathbb{Z}}(2^{ls}\|\D_{l}u\|_{L^{p}})^{q})^{\frac{1}{q}}.$$
The Besov space $B^{s}_{p,q}$ is the set of temperate distribution $u$ such that $\|u\|_{B^{s}_{p,q}}<+\infty$.
\end{definition}
\begin{remarka}The above definition is a natural generalization of the
nonhomogeneous Sobolev and H$\ddot{\mbox{o}}$lder spaces: one can show
that $B^{s}_{\infty,\infty}$ is the nonhomogeneous
H$\ddot{\mbox{o}}$lder space $C^{s}$ and that $B^{s}_{2,2}$ is
the nonhomogeneous space $H^{s}$.
\end{remarka}
\begin{proposition}
\label{derivation,interpolation}
The following properties holds:
\begin{enumerate}
\item there exists a constant universal $C$
such that:\\
$C^{-1}\|u\|_{B^{s}_{p,r}}\leq\|\n u\|_{B^{s-1}_{p,r}}\leq
C\|u\|_{B^{s}_{p,r}}.$
\item If
$p_{1}<p_{2}$ and $r_{1}\leq r_{2}$ then $B^{s}_{p_{1},r_{1}}\hookrightarrow
B^{s-N(1/p_{1}-1/p_{2})}_{p_{2},r_{2}}$.
\item $B^{s^{'}}_{p,r_{1}}\hookrightarrow B^{s}_{p,r}$ if $s^{'}> s$ or if $s=s^{'}$ and $r_{1}\leq r$.
\end{enumerate}
\label{interpolation}
\end{proposition}
Before going further into the paraproduct for Besov spaces, let us state an important proposition.
\begin{proposition}
Let $s\in\R$ and $1\leq p,r\leq+\infty$. Let $(u_{q})_{q\geq-1}$ be a sequence of functions such that
$$(\sum_{q\geq-1}2^{qsr}\|u_{q}\|_{L^{p}}^{r})^{\frac{1}{r}}<+\infty.$$
If $\mbox{supp}\hat{u}_{1}\subset {\cal C}(0,2^{q}R_{1},2^{q}R_{2})$ for some $0<R_{1}<R_{2}$ then $u=\sum_{q\geq-1}u_{q}$ belongs to $B^{s}_{p,r}$ and there exists a universal constant $C$ such that:
$$\|u\|_{B^{s}_{p,r}}\leq C^{1+|s|}\big(\sum_{q\geq-1}(2^{qs}\|u_{q}\|_{L^{p}})^{r}\big)^{\frac{1}{r}}.$$
\label{resteimp1}
\end{proposition}
Let now recall a few product laws in Besov spaces coming directly from the paradifferential calculus of J-M. Bony
(see \cite{BJM}) and rewrite on a generalized form in \cite{AP} by H. Abidi and M. Paicu (in this article the results are written
in the case of homogeneous sapces but it can easily generalize for the nonhomogeneous Besov spaces).
\begin{proposition}
\label{produit1}
We have the following laws of product:
\begin{itemize}
\item For all $s\in\R$, $(p,r)\in[1,+\infty]^{2}$ we have:
\begin{equation}
\|uv\|_{B^{s}_{p,r}}\leq
C(\|u\|_{L^{\infty}}\|v\|_{B^{s}_{p,r}}+\|v\|_{L^{\infty}}\|u\|_{B^{s}_{p,r}})\,.
\label{2.2}
\end{equation}
\item Let $(p,p_{1},p_{2},r,\lambda_{1},\lambda_{2})\in[1,+\infty]^{2}$ such that:$\frac{1}{p}\leq\frac{1}{p_{1}}+\frac{1}{p_{2}}$,
$p_{1}\leq\lambda_{2}$, $p_{2}\leq\lambda_{1}$, $\frac{1}{p}\leq\frac{1}{p_{1}}+\frac{1}{\lambda_{1}}$ and
$\frac{1}{p}\leq\frac{1}{p_{2}}+\frac{1}{\lambda_{2}}$. We have then the following inequalities:\\
if $s_{1}+s_{2}+N\inf(0,1-\frac{1}{p_{1}}-\frac{1}{p_{2}})>0$, $s_{1}+\frac{N}{\lambda_{2}}<\frac{N}{p_{1}}$ and
$s_{2}+\frac{N}{\lambda_{1}}<\frac{N}{p_{2}}$ then:
\begin{equation}
\|uv\|_{B^{s_{1}+s_{2}-N(\frac{1}{p_{1}}+\frac{1}{p_{2}}-\frac{1}{p})}_{p,r}}\lesssim\|u\|_{B^{s_{1}}_{p_{1},r}}
\|v\|_{B^{s_{2}}_{p_{2},\infty}},
\label{2.3}
\end{equation}
when $s_{1}+\frac{N}{\lambda_{2}}=\frac{N}{p_{1}}$ (resp $s_{2}+\frac{N}{\lambda_{1}}=\frac{N}{p_{2}}$) we replace
$\|u\|_{B^{s_{1}}_{p_{1},r}}\|v\|_{B^{s_{2}}_{p_{2},\infty}}$ (resp $\|v\|_{B^{s_{2}}_{p_{2},\infty}}$) by
$\|u\|_{B^{s_{1}}_{p_{1},1}}\|v\|_{B^{s_{2}}_{p_{2},r}}$ (resp $\|v\|_{B^{s_{2}}_{p_{2},\infty}\cap L^{\infty}}$),
if $s_{1}+\frac{N}{\lambda_{2}}=\frac{N}{p_{1}}$ and $s_{2}+\frac{N}{\lambda_{1}}=\frac{N}{p_{2}}$ we take $r=1$.
\\
If $s_{1}+s_{2}=0$, $s_{1}\in(\frac{N}{\lambda_{1}}-\frac{N}{p_{2}},\frac{N}{p_{1}}-\frac{N}{\lambda_{2}}]$ and
$\frac{1}{p_{1}}+\frac{1}{p_{2}}\leq 1$ then:
\begin{equation}
\|uv\|_{B^{-N(\frac{1}{p_{1}}+\frac{1}{p_{2}}-\frac{1}{p})}_{p,\infty}}\lesssim\|u\|_{B^{s_{1}}_{p_{1},1}}
\|v\|_{B^{s_{2}}_{p_{2},\infty}}.
\label{2.4}
\end{equation}
If $|s|<\NN$ for $p\geq2$ and $-\frac{N}{p^{'}}<s<\NN$ else, we have:
\begin{equation}
\|uv\|_{B^{s}_{p,r}}\leq C\|u\|_{B^{s}_{p,r}}\|v\|_{B^{\NN}_{p,\infty}\cap L^{\infty}}.
\label{2.5}
\end{equation}
\end{itemize}
\end{proposition}
\begin{remarka}
In the sequel $p$ will be either $p_{1}$ or $p_{2}$ and in this case $\frac{1}{\lambda}=\frac{1}{p_{1}}-\frac{1}{p_{2}}$
if $p_{1}\leq p_{2}$, resp $\frac{1}{\lambda}=\frac{1}{p_{2}}-\frac{1}{p_{1}}$
if $p_{2}\leq p_{1}$.
\end{remarka}
\begin{corollaire}
\label{produit2}
Let $r\in [1,+\infty]$, $1\leq p\leq p_{1}\leq +\infty$ and $s$ such that:
\begin{itemize}
\item $s\in(-\frac{N}{p_{1}},\frac{N}{p_{1}})$ if $\frac{1}{p}+\frac{1}{p_{1}}\leq 1$,
\item $s\in(-\frac{N}{p_{1}}+N(\frac{1}{p}+\frac{1}{p_{1}}-1),\frac{N}{p_{1}})$ if $\frac{1}{p}+\frac{1}{p_{1}}> 1$,
\end{itemize}
then we have if $u\in B^{s}_{p,r}$ and $v\in B^{\frac{N}{p_{1}}}_{p_{1},\infty}\cap L^{\infty}$:
$$\|uv\|_{B^{s}_{p,r}}\leq C\|u\|_{B^{s}_{p,r}}\|v\|_{B^{\frac{N}{p_{1}}}_{p_{1},\infty}\cap L^{\infty}}.$$
\end{corollaire}
The study of non stationary PDE's requires space of type $L^{\rho}(0,T,X)$ for appropriate Banach spaces $X$. In our case, we
expect $X$ to be a Besov space, so that it is natural to localize the equation through Littlewood-Payley decomposition. But, in doing so, we obtain
bounds in spaces which are not type $L^{\rho}(0,T,X)$ (except if $r=p$).
We are now going to
define the spaces of Chemin-Lerner in which we will work, which are
a refinement of the spaces
$L_{T}^{\rho}(B^{s}_{p,r})$.
$\hspace{15cm}$
\begin{definition}
Let $\rho\in[1,+\infty]$, $T\in[1,+\infty]$ and $s_{1}\in\R$. We set:
$$\|u\|_{\widetilde{L}^{\rho}_{T}(B^{s_{1}}_{p,r})}=
\big(\sum_{l\in\mathbb{Z}}2^{lrs_{1}}\|\D_{l}u(t)\|_{L^{\rho}(L^{p})}^{r}\big)^{\frac{1}{r}}\,.$$
We then define the space $\widetilde{L}^{\rho}_{T}(B^{s_{1}}_{p,r})$ as the set of temperate distribution $u$ over
$(0,T)\times\R^{N}$ such that 
$\|u\|_{\widetilde{L}^{\rho}_{T}(B^{s_{1}}_{p,r})}<+\infty$.
\end{definition}
We set $\widetilde{C}_{T}(\widetilde{B}^{s_{1}}_{p,r})=\widetilde{L}^{\infty}_{T}(\widetilde{B}^{s_{1}}_{p,r})\cap
{\cal C}([0,T],B^{s_{1}}_{p,r})$.
Let us emphasize that, according to Minkowski inequality, we have:
$$\|u\|_{\widetilde{L}^{\rho}_{T}(B^{s_{1}}_{p,r})}\leq\|u\|_{L^{\rho}_{T}(B^{s_{1}}_{p,r})}\;\;\mbox{if}\;\;r\geq\rho
,\;\;\;\|u\|_{\widetilde{L}^{\rho}_{T}(B^{s_{1}}_{p,r})}\geq\|u\|_{L^{\rho}_{T}(B^{s_{1}}_{p,r})}\;\;\mbox{if}\;\;r\leq\rho
.$$
\begin{remarka}
It is easy to generalize proposition \ref{produit1},
to $\widetilde{L}^{\rho}_{T}(B^{s_{1}}_{p,r})$ spaces. The indices $s_{1}$, $p$, $r$
behave just as in the stationary case whereas the time exponent $\rho$ behaves according to H\"older inequality.
\end{remarka}
Here we recall a result of interpolation which explains the link
of the space $B^{s}_{p,1}$ with the space $B^{s}_{p,\infty}$, see
\cite{DFourier}.
\begin{proposition}
\label{interpolationlog}
There exists a constant $C$ such that for all $s\in\R$, $\e>0$ and
$1\leq p<+\infty$,
$$\|u\|_{\widetilde{L}_{T}^{\rho}(B^{s}_{p,1})}\leq C\frac{1+\e}{\e}\|u\|_{\widetilde{L}_{T}^{\rho}(B^{s}_{p,\infty})}
\biggl(1+\log\frac{\|u\|_{\widetilde{L}_{T}^{\rho}(B^{s+\e}_{p,\infty})}}
{\|u\|_{\widetilde{L}_{T}^{\rho}(B^{s}_{p,\infty})}}\biggl).$$ \label{5Yudov}
\end{proposition}
Now we give some result on the behavior of the Besov spaces via some pseudodifferential operator (see \cite{DFourier}).
\begin{definition}
Let $m\in\R$. A smooth function function $f:\R^{N}\rightarrow\R$ is said to be a ${\cal S}^{m}$ multiplier if for all muti-index $\alpha$, there exists a constant $C_{\alpha}$ such that:
$$\forall\xi\in\R^{N},\;\;|\p^{\alpha}f(\xi)|\leq C_{\alpha}(1+|\xi|)^{m-|\alpha|}.$$
\label{smoothf}
\end{definition}
\begin{proposition}
Let $m\in\R$ and $f$ be a ${\cal S}^{m}$ multiplier. Then for all $s\in\R$ and $1\leq p,r\leq+\infty$ the operator $f(D)$ is continuous from $B^{s}_{p,r}$ to $B^{s-m}_{p,r}$.
\label{singuliere}
\end{proposition}
\section{Estimates for parabolic system with variable coefficients}
\label{section3}
Let us first state estimates for the following constant coefficient parabolic system:
\begin{equation}
\begin{cases}
\begin{aligned}
&\p_{t}u-\mu\D u-(\lambda+\mu)\n{\rm div}u=f,\\
&u_{/t=0}=u_{0}.
\end{aligned}
\end{cases}
\label{3}
\end{equation}
\begin{proposition}
Assume that $\mu\geq0$ and that $\lambda+2\mu\geq0$. Then there exists a universal constant $\kappa$ such that for all $s\in\mathbb{Z}$ and $T\in\R^{+}$,
$$
\begin{aligned}
&\|u\|_{\widetilde{L}^{\infty}_{T}(B^{s}_{p_{1},1})}\leq\|u_{0}\|_{B^{s}_{p_{1},1}}+\|f\|_{L^{1}_{T}(B^{s}_{p_{1},1})},\\
&\kappa\nu\|u\|_{L^{1}_{T}(B^{s+2}_{p_{1},1})}\leq\sum_{l\in\mathbb{Z}}2^{ls}(1-e^{-\kappa\nu2^{2l}T})(\|\D_{l}u_{0}\|_{L^{p_{1}}}
+\|\D_{l}f\|_{L^{1}_{T}(L^{p_{1}})}),
\end{aligned}
$$
with $\nu=\min(\mu,\lambda+2\mu)$.
\label{chaleur}
\end{proposition}
We now consider the following parabolic system which is obtained by linearizing the momentum equation:
\begin{equation}
\begin{cases}
\begin{aligned}
&\p_{t}u+v\cdot\n u+u\cdot\n w-b(\mu\D u+(\lambda+\mu)\n{\rm div}u=f+g,\\
&u_{/t=0}=u_{0}.
\label{5}
\end{aligned}
\end{cases}
\end{equation}
Above $u$ is the unknown function. We assume that $u_{0}\in B^{s}_{p_{1},1}$, $f\in L^{1}(0,T;B^{s}_{p_{1},1})$ and $g\in L^{r}(0,T;B^{s^{'}}_{q_{1},1})$, that $v$ and $w$ are time dependent
vector-fields with coefficients in $L^{1}(0,T;B^{\NN+1}_{p,1})$, that $b$ is bounded by below by a positive constant $\underline{b}$ and 
$b$ belongs to $L^{\infty}(0,T;B^{\frac{N}{p}}_{p,1})$ with $p\in[1,+\infty]$.
\begin{proposition}
Let $g=0$ and $\underline{\nu}=\underline{b}\min(\mu,\lambda+2\mu)$ and $\bar{\nu}=\mu+|\lambda+\mu|$. Assume that $s\in(-\frac{N}{p},\frac{N}{p}]$. Let $m\in\mathbb{Z}$
be such that $b_{m}=1+S_{m}a$ satisfies:
\begin{equation}
\inf_{(t,x)\in[0,T)\times\R^{N}}b_{m}(t,x)\geq\frac{\underline{b}}{2}.
\label{6}
\end{equation}
There exist three constants $c$, $C$ and $\kappa$ (with $c$, $C$, depending only on $N$ and on $s$, and $\kappa$ universal) such that if in addition we have:
\begin{equation}
\|1-S_{m}a\|_{L^{\infty}(0,T;B^{\frac{N}{p_{1}}}_{p_{1},1})}\leq c\frac{\underline{\nu}}{\bar{\nu}}
\label{7}
\end{equation}
then setting:
$$V(t)=\int^{t}_{0}\|v\|_{B^{\NN+1}_{p,1}}d\tau,\;\;W(t)=\int^{t}_{0}\|w\|_{B^{\NN+1}_{p,1}}d\tau,\;\;\mbox{and}
\;\;Z_{m}(t)=2^{2m}\bar{\nu}^{2}\underline{\nu}^{-1}\int^{t}_{0}\|a\|^{2}_{B^{\frac{N}{p_{1}}}_{p_{1},1}}d\tau,$$
We have for all $t\in[0,T]$,
$$
\begin{aligned}
&\|u\|_{\widetilde{L}^{\infty}((0,T)\times B^{s}_{p_{1},1})}+\kappa\underline{\nu}
\|u\|_{\widetilde{L}^{1}((0,T)\times B^{s+2}_{p_{1},1})}\leq e^{C(V+W+Z_{m})(t)}(\|u_{0}\|_{B^{s}_{p_{1},1}}\\
&\hspace{7cm}+\int^{t}_{0}
e^{-C(V+W+Z_{m})(\tau)}\|f(\tau)\|_{B^{s}_{p_{1},1}}d\tau).
\end{aligned}
$$
\label{linearise}
\end{proposition}
\begin{remarka}
Let us stress the fact that if $a\in \widetilde{L}^{\infty}((0,T)\times B^{\NN}_{p,1})$ then assumption (\ref{6}) and
(\ref{7}) are satisfied for $m$ large enough. This will be used in the proof of theorem \ref{theo1}.
Indeed, according to Bernstein inequality, we have:
$$\|a-S_{m}a\|_{L^{\infty}((0,T)\times\R^{N})}\leq\sum_{q\geq m}\|\D_{q}a\|_{L^{\infty}((0,T)\times\R^{N})}\lesssim\sum_{q\geq m}
2^{q\NN}\|\D_{q}a\|_{L^{\infty}(L^{p})}.$$
Because $a\in\widetilde{L}^{\infty}((0,T)\times B^{\NN}_{p,1})$, the right-hand side is the remainder of a convergent series hence tends to zero
when $m$ goes to infinity. For a similar reason, (\ref{7}) is satisfied for $m$ large enough.
\label{remark5}
\end{remarka}
{\bf Proof:}
Let us first rewrite (\ref{5}) as follows:
\begin{equation}
\p_{t}u+v\cdot\n u+u\cdot\n w-b_{m}(\mu\D u+(\lambda+\mu)\n{\rm div}u=f+E_{m}-u\cdot\n w,
\label{8}
\end{equation}
with $E_{m}=(\mu\D u+(\lambda+\mu)\n{\rm div}u)(\mbox{Id}-S_{m})a$.
Note that, because $-\NN<s\leq\NN$, the error term $E_{m}$ may be estimated by:
\begin{equation}
\|E_{m}\|_{B^{s}_{p_{1},1}}\lesssim\|a-S_{m}a\|_{B^{\NN}_{p,1}}\|D^{2}u\|_{B^{s}_{p_{1},1}}.
\label{9}
\end{equation}
and we have:
\begin{equation}
\|u\cdot\n w\|_{B^{s}_{p_{1},1}}\lesssim\|\n w\|_{B^{\NN}_{p,1}}\|u\|_{B^{s}_{p_{1},1}}.
\label{10}
\end{equation}
Now applying $\D_{q}$ to equation (\ref{8}) yields:
\begin{equation}
\begin{aligned}
\frac{d}{dt}u_{q}+v\cdot\n u_{q}-\mu{\rm div}(b_{m}\n u_{q})-(\lambda+\mu)\n(b_{m}{\rm div}u_{q})=f_{q}&+E_{m,q}-\D_{q}(u\cdot\n w)\\
&\hspace{1,5cm}+R_{q}+\widetilde{R}_{q}.
\end{aligned}
\end{equation}
where we denote by $u_{q}=\D_{q}u$ and
with:
$$
\begin{aligned}
&R_{q}=[v^{j},\D_{q}]\p_{j}u,\\
&\widetilde{R}_{q}=\mu\big(\D_{q}(b_{m}\D u)-{\rm div}(b_{m}\n u_{q})\big)+(\lambda+\mu)\big(\D_{q}(b_{m}\n{\rm div}u)-\n(b_{m}{\rm div}u_{q})\big).
\end{aligned}
$$
Next multiplying both sides by $|u_{q}|^{p_{1}-2}u_{q}$, and integrating by parts in the second, third and last term in the left-hand side, we get:
$$
\begin{aligned}
&\frac{1}{p_{1}}\frac{d}{dt}\|u_{q}\|_{L^{p_{1}}}^{p_{1}}-\frac{1}{p_{1}}\int\big(|u_{q}|^{p_{1}}{\rm div}v+\mu
{\rm div}(b_{m}\n u_{q})|u_{q}|^{p_{1}-2}u_{q}
+\xi\n\big(b_{m}{\rm div}u_{q}\big)|u_{q}|^{p_{1}-2}u_{q})\big)dx\\
&\hspace{1,7cm}\leq\|u_{q}\|^{p_{1}-1}_{L^{p_{1}}}(\|f_{q}\|_{L^{p_{1}}}+\|\D_{q}E_{m}\|_{L^{p_{1}}}+\|\D_{q}(u\cdot\n w)\|_{L^{p_{1}}}+\|R_{q}\|_{L^{p_{1}}}
+\|\widetilde{R}_{q}\|_{L^{p_{1}}}).
\end{aligned}
$$
Hence denoting $\xi=\mu+\lambda$, $\nu=\min(\mu,\lambda+2\mu)$ and using (\ref{6}), lemma [A5] of \cite{DL} and Young's inequalities we get:
$$
\begin{aligned}
&\frac{1}{p_{1}}\frac{d}{dt}\|u_{q}\|_{L^{p_{1}}}^{p_{1}}+\frac{\nu \underline{b}(p_{1}-1)}{p_{1}^{2}}2^{2q}\|u_{q}\|_{L^{p_{1}}}^{p_{1}}\leq
\|u_{q}\|^{p_{1}-1}_{L^{p_{1}}}\big(\|f_{q}\|_{L^{p_{1}}}+\|E_{m,q}\|_{L^{p_{1}}}
+\|\D_{q}(u\cdot\n w)\|_{L^{p_{1}}}\\
&\hspace{6,8cm}+\frac{1}{p_{1}}\|u_{q}\|_{L^{p_{1}}}\|{\rm div}u\|_{L^{\infty}}+\|R_{q}\|_{L^{p_{1}}}+\|\widetilde{R}_{q}\|_{L^{p_{1}}}\big),
\end{aligned}
$$
which leads, after time integration to:
\begin{equation}
\begin{aligned}
&\|u_{q}\|_{L^{p_{1}}}+\frac{\nu \underline{b}(p_{1}-1)}{p_{1}}2^{2q}\int^{t}_{0}\|u_{q}\|_{L^{p_{1}}}d\tau\leq
\|\D_{q}u_{0}\|_{L^{p_{1}}}+\int^{t}_{0}\big(
\|f_{q}\|_{L^{p_{1}}}+\|E_{m,q}\|_{L^{p_{1}}}\\
&\hspace{1,9cm}+\|\D_{q}(u\cdot\n w)\|_{L^{p_{1}}}+\frac{1}{p_{1}}\|u_{q}\|_{L^{p_{1}}}\|{\rm div}u\|_{L^{\infty}}+\|R_{q}\|_{L^{p_{1}}}+
\|\widetilde{R}_{q}\|_{L^{p_{1}}}\big)d\tau,
\end{aligned}
\label{11}
\end{equation}
where $\underline{\nu}=\underline{b}\nu$.
For commutators $R_{q}$ and $\widetilde{R}_{q}$, we have the following estimates (see lemma \ref{alemme2} and \ref{alemme3}
in the appendix)
\begin{equation}
\|R_{q}\|_{L^{p_{1}}}\lesssim c_{q}2^{-qs}\|v\|_{B^{\frac{N}{p}+1}_{p,1}}\|u\|_{B^{s}_{p_{1},1}},\label{12}
\end{equation}
\begin{equation}
\|\widetilde{R}_{q}\|_{L^{p_{1}}}\lesssim c_{q}\bar{\nu}2^{-qs}\|S_{m}a\|_{B^{\frac{N}{p_{1}}+1}_{p_{1},1}}\|Du\|_{B^{s}_{p_{1},1}},\label{13}
\end{equation}
where $(c_{q})_{q\in\mathbb{Z}}$ is a positive sequence such that $\sum_{q\in\mathbb{Z}}c_{q}=1$, and $\bar{\nu}=\mu+|\lambda+\mu|$. Note that,
owing to Bernstein inequality, we have:
$$\|S_{m}a\|_{B^{\NN+1}_{p,1}}\lesssim2^{m}\|a\|_{B^{\frac{N}{p_{1}}}_{p_{1},1}}$$
Hence, plugging these latter estimates and (\ref{9}), (\ref{10}) in (\ref{11}), then multiplying by $2^{qs}$ and summing up on $q\in\mathbb{Z}$, we discover that, for all $t\in[0,T]$:
$$
\begin{aligned}
&\|u\|_{L^{\infty}_{t}(B^{s}_{p_{1},1})}+\frac{\nu \underline{b}(p_{1}-1)}{p}\|u\|_{L^{1}_{t}(B^{s+2}_{p_{1},1})}\leq
\|u_{0}\|_{B^{s}_{p_{1},1}}+
\|f\|_{L^{1}_{t}(B^{s}_{p_{1},1})}+C\int^{t}_{0}(\|v\|_{B^{\frac{N}{p}+1}_{p_{1},1}}\\
&\hspace{1cm}+\|w\|_{B^{\frac{N}{p}+1}_{p,1}})\|u\|_{B^{s}_{p_{1},1}}d\tau+C\bar{\nu}\int^{t}_{0}(\|a-S_{m}a\|_{B^{\NN}_{p,1}}
\|u\|_{B^{s+2}_{p_{1},1}}
+2^{m}\|a\|_{B^{\NN}_{p,1}}\|u\|_{B^{s+1}_{p_{1},1}})d\tau,
\end{aligned}
$$
for a constant $C$ depending only on $N$ and $s$.
Let $X(t)=\|u\|_{L^{\infty}_{t}(B^{s}_{p_{1},1})}+\nu \underline{b}\|u\|_{L^{1}_{t}(B^{s+2}_{p_{1},1})}$. Assuming that $m$ has been chosen so
large as to satisfy:
$$C\bar{\nu}\|a-S_{m}a\|_{L^{\infty}_{T}(B^{\NN}_{p,1})}\leq\underline{\nu},$$
and using that by interpolation, we have:
$$C\bar{\nu}\|a\|_{B^{\NN}_{p,1}}\|u\|_{B^{s+2}_{p,1}}\leq\kappa\underline{\nu}+\frac{C^{2}\bar{\nu}^{2}2^{2m}}
{4\kappa\underline{\nu}}
\|a\|^{2}_{B^{\NN}_{p,1}}\|u\|_{B^{s}_{p,1}},$$
we end up with:
$$X(t)\leq\|u_{0}\|_{B^{s}_{p_{1},1}}+
\|f\|_{L^{1}_{t}(B^{s}_{p_{1},1})}+C\int^{t}_{0}(\|v\|_{B^{\NN+1}_{p,1}}+\|w\|_{B^{\NN+1}_{p,1}}+\frac{\bar{\nu}^{2}}
{\underline{\nu}}2^{2m}
\|a\|^{2}_{B^{\N}_{p,1}})Xd\tau\\$$
Gr\"onwall lemma then leads to the desired inequality.
{\hfill $\Box$}
\begin{remarka}
The proof of the continuation criterion (theorem \ref{theo1}) relies on a better estimate which is available when $u=v=w$ and $s>0$. In fact, by arguing as in the
proof of the previous proposition and by making use of inequality (\ref{54}) instead of (\ref{52}), one can prove that
under conditions (\ref{6}) and (\ref{7}),
there exists constants $C$ and $\kappa$ such that:
$$
\begin{aligned}
&\forall t\in[0,T],\;\|u\|_{L^{\infty}_{t}(B^{s}_{p_{1},1})}+\kappa\underline{\nu}\|u\|_{L^{1}_{t}(B^{s+2}_{p_{1},1})}
\leq e^{C(U+Z_{m})(t)}
\big(\|u_{0}\|_{B^{s}_{p_{1},1}}+\\
&\hspace{3cm}\int^{t}_{0}e^{-C(U+Z_{m})(\tau)}\|f(\tau)\|_{B^{s}_{p_{1},1}}d\tau\big)\;\;\;\mbox{with}\;\;\;U(t)=\int^{t}_{0}\|\n u\|_{L^{\infty}}d\tau.
\end{aligned}
$$
\label{remark6}
\end{remarka}
In the following corollary, we generalize proposition \ref{linearise1} when $g\ne 0$ and $g\in \widetilde{L}^{r}(B^{s^{'}}_{q_{1},1})$.
Moreover here $u_{0}=u_{1}+u_{2}$ with $u_{1}\in B^{s}_{p_{1},1}$ and $u_{2}\in B^{s^{'}}_{p_{2},1}$.
\begin{corollaire}
Let $\underline{\nu}=\underline{b}\min(\mu,\lambda+2\mu)$ and $\bar{\nu}=\mu+|\lambda+\mu|$. Assume that $s,s^{'}\in(-\frac{N}{p},\frac{N}{p}]$. Let $m\in\mathbb{Z}$
be such that $b_{m}=1+S_{m}a$ satisfies:
\begin{equation}
\inf_{(t,x)\in[0,T)\times\R^{N}}b_{m}(t,x)\geq\frac{\underline{b}}{2}.
\label{6}
\end{equation}
There exist three constants $c$, $C$ and $\kappa$ (with $c$, $C$, depending only on $N$ and on $s$, and $\kappa$ universal) such that if in addition we have:
\begin{equation}
\|1-S_{m}a\|_{L^{\infty}(0,T;B^{\frac{N}{p}}_{p,1})}\leq c\frac{\underline{\nu}}{\bar{\nu}}
\label{7}
\end{equation}
then setting:
$$V(t)=\int^{t}_{0}\|v\|_{B^{\frac{N}{p}+1}_{p,1}}d\tau,\;\;\;W(t)=\int^{t}_{0}\|w\|_{B^{\frac{N}{p}+1}_{p,1}}d\tau,\;\;\;\mbox{and}
\;\;\;Z_{m}(t)=2^{2m}\bar{\nu}^{2}\underline{\nu}^{-1}\int^{t}_{0}\|a\|^{2}_{B^{\frac{N}{p}}_{p,1}}d\tau,$$
We have for all $t\in[0,T]$,
$$
\begin{aligned}
&\|u\|_{\widetilde{L}^{\infty}_{T}(B^{s}_{p_{1},1}+B^{s^{'}}_{p_{2},1})}+\kappa\underline{\nu}
\|u\|_{\widetilde{L}^{1}_{T}(B^{s+2}_{p_{1},1}+B^{s^{'}+2}_{p_{2},1})}\leq e^{C(V+W+Z_{m})(t)}\big(\|u_{1}\|_{B^{s}_{p_{1},1}}+\\
&\hspace{4cm}\|u_{2}\|_{B^{s^{'}}_{p_{2},1}}+\int^{t}_{0}
e^{-C(V+W+Z_{m})(\tau)}(\|f(\tau)\|_{B^{s}_{p_{1},1}}+\|g(\tau)\|_{B^{s^{'}}_{q_{1},1}})d\tau\big).
\end{aligned}
$$
\label{coroimportant}
\end{corollaire}
{\bf Proof:} We split the solution $u$ in two parts $u_{1}$ and $u_{2}$ which verify the following equations:
$$
\begin{cases}
\begin{aligned}
&\p_{t}u_{1}+v\cdot\n u_{1}+u_{1}\cdot\n w-b(\mu\D u_{1}+(\lambda+\mu)\n{\rm div}u_{1}=f,\\
&u_{/t=0}=u_{1}^{0},
\label{5}
\end{aligned}
\end{cases}
$$
and:
$$
\begin{cases}
\begin{aligned}
&\p_{t}u_{2}+v\cdot\n u_{2}+u_{2}\cdot\n w-b(\mu\D u_{2}+(\lambda+\mu)\n{\rm div}u_{2}=g,\\
&u_{/t=0}=u_{2}^{0}.
\label{5}
\end{aligned}
\end{cases}
$$
We have then $u=u_{1}+u_{2}$ and we conclude by applying proposition \ref{linearise}.
{\hfill $\Box$}\\
Proposition \ref{linearise} fails in the limit case $s=-\NN$. The reason why is that proposition \ref{produit1} cannot
be applied any longer. One can however state the following
result which will be the key to the proof of uniqueness in dimension two.
\begin{proposition}
\label{linearise1}
Under condition (\ref{6}), there exists three constants $c$, $C$ and $\kappa$ (with $c$, $C$, depending
only on $N$, and $\kappa$ universal) such that if:
\begin{equation}
\|a-S_{m}a\|_{\widetilde{L}^{\infty}_{t}(B^{\NN}_{p,1})}\leq c\frac{\underline{\nu}}{\bar{\nu}},
\label{14}
\end{equation}
then we have:
$$\|u\|_{L^{\infty}_{t}(B^{-\frac{N}{p_{1}}}_{p_{1},\infty})}+\kappa\underline{\nu}\|u\|_{\widetilde{L}^{1}_{t}(B^{2-\frac{N}{p_{1}}}
_{p_{1},\infty})}
\leq
2e^{C(V+W)(t)}(\|u_{0}\|_{B^{-\frac{N}{p_{1}}}_{p_{1},\infty}}+\|f\|_{\widetilde{L}^{1}_{t}(B^{\frac{N}{p_{1}}}_{p_{1},\infty})}),$$
whenever $t\in[0,T]$ satisfies:
\begin{equation}
\bar{\nu}^{2}t\|a\|^{2}_{\widetilde{L}^{\infty}_{t}(B^{\NN}_{p,1})}\leq c2^{-2m}\underline{\nu}.
\label{15}
\end{equation}
\end{proposition}
{\bf Proof:}
We just point out the changes that have to be be done compare to the proof of proposition \ref{linearise}. The first one is that instead
of (\ref{9}) and (\ref{10}), we have in accordance with proposition \ref{produit1}:
\begin{equation}
\|E_{m}\|_{\widetilde{L}^{1}_{t}(B^{-\frac{N}{p_{1}}}_{p_{1},\infty})}\lesssim\|a-S_{m}a\|
_{\widetilde{L}^{\infty}_{t}(B^{\NN}_{p,1})}\|D^{2}u\|_{\widetilde{L}^{1}_{t}(B^{-\frac{N}{p_{1}}}_{p_{1},\infty})},
\label{16}
\end{equation}
\begin{equation}
\|u\cdot w\|_{B^{-\NN}_{p,\infty}}\lesssim\|u\|_{B^{-\frac{N}{p_{1}}}_{p_{1},\infty}}\|\n w\|_{B^{\NN}_{p,1}}.
\label{17}
\end{equation}
The second change concerns the estimates of commutator $R_{q}$ and $\widetilde{R}_{q}$. According to inequality (\ref{53})
and remark \ref{remarque7},
we now have for all $q\in\mathbb{Z}$:
\begin{equation}
\|R_{q}\|_{L^{p}}\lesssim 2^{q\frac{N}{p_{1}}}\|v\|_{B^{\NN+1}_{p,1}}\|u\|_{B^{-\frac{N}{p_{1}}}_{p_{1},\infty}},
\label{18}
\end{equation}
\begin{equation}
\|\widetilde{R}_{q}\|\lesssim\bar{\nu}2^{q\frac{N}{p_{1}}}\|S_{m}a\|_{\widetilde{L}^{\infty}_
{t}(B^{\NN+1}_{p,1})}\|Du\|_{\widetilde{L}^{1}_{t}(B^{-\frac{N}{p_{1}}}_{p_{1},\infty})}.
\label{19}
\end{equation}
Plugging all these estimates in (\ref{11}) then taking the supremum over $q\in\mathbb{Z}$, we get:
$$
\begin{aligned}
&\|u\|_{L^{\infty}_{t}(B^{-\frac{N}{p_{1}}}_{p_{1},\infty})}+2\underline{\nu}\|u\|_{\widetilde{L}^{1}_{t}(B^{2-\frac{N}{p_{1}}}_{p_{1},
\infty})}\leq
\|u_{0}\|_{B^{-\frac{N}{p_{1}}}_{p_{1},1}}+C_int^{t}_{0}(\|v\|_{B^{\NN+1}_{p,1}}+\|w\|_{B^{\NN+1}_{p,1}})\|u\|_{B^{-\frac{N}{p_{1}}}_{p_{1},\infty}}d\tau\\
&+C\bar{\nu}\big(\|a-S_{m}a\|_{\widetilde{L}^{\infty}_{t}(B^{\NN}_{p,1})}\|u\|_{\widetilde{L}^{1}_{t}(B^{2-\frac{N}{p_{1}}}_{p_{1}
,\infty})}
+2^{m}\|a\|_{L^{\infty}_{t}(B^{\NN}_{p,1})}\|u\|_{\widetilde{L}^{1}_{t}(B^{1-\frac{N}{p_{1}}}_{p_{1},\infty})}+
\|f\|_{\widetilde{L}^{1}_{t}(B^{-\frac{N}{p_{1}}}_{p_{1},\infty})}.
\end{aligned}
$$
Using that:
$$\|u\|_{\widetilde{L}^{1}_{t}(B^{1-\frac{N}{p_{1}}}_{p_{1},\infty})}\leq\sqrt{t}\|u\|^{\frac{1}{2}}_{\widetilde{L}^{1}_{t}(B^{2-
\frac{N}{p_{1}}}_{p_{1},\infty})}\big)
\|u\|^{\frac{1}{2}}_{L^{\infty}_{t}(B^{\frac{N}{p_{1}}}_{p_{1},\infty})},$$
and taking advantage of assumption (\ref{14}) and (\ref{15}), it is now easy to complete the proof.
{\hfill $\Box$}
\section{The mass conservation equation}
\label{section4}
Let us first recall standard estimates in Besov spaces for the following linear transport equation:
$$
\begin{cases}
\begin{aligned}
&\p_{t}a+u\cdot\n a=g,\\
&a_{/t=0}=a_{0}.
\end{aligned}
\end{cases}
\leqno{({\cal H})}
$$
\begin{proposition}
Let $1\leq p_{1}\leq p\leq+\infty$, $r\in[1,+\infty]$ and $s\in\R$ be such that:
$$-N\min(\frac{1}{p_{1}},\frac{1}{p^{'}})<s<1+\frac{N}{p_{1}}.$$
There exists a constant $C$ depending only on $N$, $p$, $p_{1}$, $r$ and $s$ such that for all $a\in L^{\infty}([0,T],B^{\sigma}_{p,r})$ of $({\cal H})$ with initial data $a_{0}$ in $B^{s}_{p,r}$ and $g\in L^{1}([0,T], B^{s}_{p,r})$, we have for a.e $t\in[0,T]$:
\begin{equation}
\|f\|_{\widetilde{L}^{\infty}_{t}(B^{s}_{p,r})}\leq e^{CU(t)}\big(\|f_{0}\|_{B^{s}_{p,r}}+\int^{t}_{0}e^{-CV(\tau)}
\|F(\tau)\|_{B^{s}_{p_{1},r}}d\tau\big),
\label{20}
\end{equation}
with:
$U(t)=\int^{t}_{0}\|\n u(\tau)\|_{B^{\frac{N}{p_{1}}}_{p_{1},\infty}\cap L^{\infty}}d\tau$.
\label{transport1}
\end{proposition}
For the proof of proposition \ref{transport1}, see \cite{BCD}.
We now focus on the mass equation associated to (\ref{0.6}):
\begin{equation}
\begin{cases}
\begin{aligned}
&\p_{t}a+v\cdot\n a=(1+a){\rm div}v,\\
&a_{/t=0}=a_{0}.
\end{aligned}
\end{cases}
\label{525}
\end{equation}
Here we generalize a proof of R. Danchin in \cite{DW}.
\begin{proposition}
Let  $r\in{1,+\infty}$, $1\leq p_{1}\leq p\leq+\infty$ and $s\in(-\min(\frac{N}{p_{1}},\frac{N}{p^{'}},\NN]$. Assume that $a_{0}\in B^{s}_{p,r}\cap L^{\infty}$, $v\in L^{1}(0,T;B^{\frac{N}{p_{1}}+1}_{p_{1},1})$ and that
$a\in\widetilde{L}^{\infty}_{T}(B^{s}_{p,r})\cap L^{\infty}_{T}$ satisfies (\ref{525}).
Let $V(t)=\int^{t}_{0}\|\n v(\tau)\|_{B^{\frac{N}{p_{1}}}_{p_{1},1}}d\tau$.
There exists a constant $C$ depending only on $N$ such that for all
$t\in[0,T]$ and $m\in\mathbb{Z}$, we have:
\begin{equation}
\|a\|_{\widetilde{L}^{\infty}_{t}(B^{s}_{p,r}\cap L^{\infty})}\leq e^{2CV(t)}\|a_{0}\|_{B^{s}_{p,r}\cap L^{\infty}}+e^{2CV(t)}-1,
\label{22}
\end{equation}
\begin{equation}
\|a-S_{m}a\|_{B^{s}_{p,r}}\leq\|a_{0}-S_{m}a_{0}\|_{B^{s}_{p,r}}+
\frac{1}{2}(1+\|a_{0}\|_{B^{s}_{p,r}\cap L^{\infty}})(e^{2CV(t)}-1)+C\|a\|_{L^{\infty}}V(t),
\label{23}
\end{equation}
\begin{equation}
\begin{aligned}
&\big(\sum_{l\leq m}2^{lrs}\|\D_{l}(a-a_{0})\|^{r}_{L^{\infty}_{t}(L^{p})}\big)^{\frac{1}{r}}\leq(1+\|a_{0}\|_{B^{s}_{p,r}})(e^{CV(t)}-1)
\\
&\hspace{8cm}+C2^{m}\|a_{0}\|_{B^{s}_{p,r}}\int^{t}_{0}\|v\|_{B^{\frac{N}{p_{1}}}_{p_{1},1}}d\tau.
\label{24}
\end{aligned}
\end{equation}
\label{transport2}
\end{proposition}
{\bf Proof:}\;\;
Applying $\D_{l}$ to (\ref{525}) yields:
$$\p_{t}\D_{l}a+v\cdot\n\D_{l}a=R_{l}+\D_{l}((1+a){\rm div}v)\;\;\;\mbox{with}\;\;R_{l}=[v\cdot\n,\D_{l}]a.$$
Multipling by $\D_{l}a|\D_{l}a|^{p-2}$ then performing a time integration, we easily get:
$$\|\D_{l}a(t)\|_{L^{p}}\lesssim\|\D_{l}a_{0}\|_{L^{p}}+\int^{t}_{0}\big(\|R_{l}\|_{L^{p}}+\|{\rm div}v\|_{L^{\infty}}\|\D_{l}a\|_{L^{p}}
+\|\D_{l}((1+a){\rm div}v)\|_{L^{p}}\big)d\tau.$$
According to proposition \ref{produit1} and interpolation, there exists a constant $C$ and a positive sequence $(c_{l})_{l\in\mathbb{N}}$ in $l^{r}$ with norm $1$ such that:
$$\|\D_{l}((1+a){\rm div}v)\|_{L^{p}}\leq Cc_{l}2^{-ls}(1+\|a\|_{B^{s}_{p,r}\cap L^{\infty}})\|{\rm div}v\|_{B^{\frac{N}{p_{1}}}_{p_{1},1}}.$$
Next the term $\|R_{l}\|_{L^{p}}$ may be bounded according to lemma \ref{alemme2} in appendix. We end up with:
\begin{equation}
\forall t\in [0,T],\;\forall l\in\mathbb{Z},\;\;2^{ls}\|\D_{l}a(t)\|_{L^{p}}\leq 2^{ls}\|\D_{l}a_{0}\|_{L^{p}}+
C\int^{t}_{0}c_{l}(1+\|a\|_{B^{s}_{p,r}\cap L^{\infty}})V^{'}d\tau,
\label{25}
\end{equation}
hence, summing up on $\mathbb{Z}$ in $l^{r}$,
$$\forall t\in [0,T],\;\forall l\in\mathbb{Z},\;\;\|a(t)\|_{B^{s}_{p,r}}\leq\|a_{0}\|_{B^{s}_{p,r}}+\int^{t}_{0}CV^{'}\|a(\tau)\|_{B^{s}_{p,r}}d\tau
+\int^{t}_{0}C(1+\|a\|_{L^{\infty}_{T}})V^{'}d\tau.
$$
Next we have:
$$\|a\|_{L^{\infty}_{t}}\leq\int^{t}_{0}(1+\|a(\tau)\|_{L^{\infty}})V^{'}(\tau)d\tau.$$
By summing the two previous inequalities, applying Gronwall lemma and proposition \ref{resteimp1} yields inequality (\ref{22}).
Let us now prove inequality (\ref{23}). Starting from (\ref{25}) and summing up over $l\geq m$ in $l^{r}$, we get:
$$
\begin{aligned}
&(\sum_{l\geq m}2^{lsr}\|\D_{l}a\|^{r}_{L^{\infty}_{t}(L^{p})})^{\frac{1}{r}}\leq (\sum_{l\geq m}2^{lsr}\|\D_{l}a_{0}\|^{r}_{L^{p}})^{\frac{1}{r}}+
C\int^{t}_{0}V^{'}(e^{2CV}\|a_{0}\|_{B^{s}_{p,r}\cap L^{\infty}}+e^{2CV}-1)d\tau\\
&\hspace{10cm}+\int^{t}_{0}C(1+\|a\|_{L^{\infty}})V^{'}d\tau.
\end{aligned}
$$
Straightforward calculations then leads to (\ref{23}).
In order to prove (\ref{24}), we use the fact that $\widetilde{a}=a-a_{0}$ satisfies:
$$\p_{t}\widetilde{a}+v\cdot\n\widetilde{a}=(1+\widetilde{a}){\rm div}v+a_{0}{\rm div}v-v\cdot\n a_{0},\;\;\widetilde{a}_{/t=0}=0.$$
Therefore, arguing as for proving (\ref{25}), we get for all $t\in[0,T]$ and $l\in\mathbb{Z}$,
$$
\begin{aligned}
&2^{l\NN}\|\D_{l}\widetilde{a}\|_{L^{p}}\leq \int^{t}_{0}2^{l\NN}\big(\|\D_{l}(a_{0}{\rm div}v)\|_{L^{p}}+\|\D_{l}(v\cdot\n a_{0})\|_{L^{p}}\big)d\tau\\
&\hspace{8cm}+C\int^{t}_{0}c_{l}(1+\|a\|_{B^{\NN}_{p,1}})V^{'}d\tau.
\end{aligned}
$$
Since $B^{\NN}_{p,1}$ is an algebra and the product maps $B^{\NN}_{p,1}\times B^{\NN-1}_{p,1}$ in $B^{\NN-1}_{p,1}$, we discover that:
$$
\begin{aligned}
&2^{l\NN}\|\D_{l}\widetilde{a}\|_{L^{\infty}(L^{p})}\leq C\big(\int^{t}_{0}2^{l}c_{l}\|a_{0}\|_{B^{\NN}_{p,1}}\|v\|_{B^{\NN}_{p,1}}d\tau+
\int^{t}_{0}c_{l}(1+\|a_{0}\|_{B^{\NN}_{p,1}}+\|a\|_{B^{\NN}_{p,1}})V^{'}d\tau\big),
\end{aligned}
$$
hence, summing up on $l\leq m$,
$$\begin{aligned}
&\sum_{l\leq m}2^{l\NN}\|\D_{l}\widetilde{a}\|_{L^{\infty}(L^{p})}\leq C\big(\int^{t}_{0}2^{m}\|a_{0}\|_{B^{\NN}_{p,1}}\|v\|_{B^{\NN}_{p,1}}d\tau+
\int^{t}_{0}(1+\|a_{0}\|_{B^{\NN}_{p,1}}+\|a\|_{B^{\NN}_{p,1}})V^{'}d\tau\big),
\end{aligned}
$$
Plugging (\ref{22}) in the right-hand side yields (\ref{24}).
\section{The proof of theorem \ref{theo1}}
\label{section5}
\subsection{Strategy of the proof}
To improve the results of R. Danchin in \cite{DL}, \cite{DW}, it is crucial to kill the coupling between the velocity and the pressure which intervene in the works of R. Danchin. In this goal, we need to integrate the pressure term in the study of the linearized equation of the momentum equation. For making, we will try to express the gradient of the pressure as a Laplacian term, so we set for $\bar{\rho}>0$ a constant state:
$${\rm div}v=P(\rho)-P(\bar{\rho}).$$
Let ${\cal E}$ the fundamental solution of the Laplace operator.
$$$$
We will set in the sequel: $v=\n{\cal E}*\big(P(\rho)-P(\bar{\rho})\big)=\n\big({\cal E}*[P(\rho)-P(\bar{\rho})]\big)$ ( $*$ here means the operator of convolution). We verify next that:
$$
\begin{aligned}
\n{\rm div}v=\n\D \big({\cal E}*[P(\rho)-P(\bar{\rho})]\big)=\D\n\big({\cal E}*[P(\rho)-P(\bar{\rho})]\big)=\D v=\n P(\rho).
\end{aligned}
$$
By this way we can now rewrite the momentum equation of (\ref{0.6}). We obtain the following equation where we have set $\nu=2\mu+\lambda$:
$$\p_{t}u+u\cdot \n u-\frac{\mu}{\rho}\D\big(u-\frac{1}{\nu}v\big)-\frac{\lambda+\mu}{\rho}\n{\rm div}\big(u-\frac{1}{\nu}v\big)=f.$$
We want now calculate $\p_{t}v$, by the transport equation we get:
$$\p_{t}v=\n{\cal E}*\p_{t}P(\rho)=-\n {\cal E}*\big(P^{'}(\rho){\rm div}(\rho u)\big).$$
We have finally:
$$\D(\p_{t}F)=-P^{'}(\rho){\rm div}(\rho u).$$
\begin{notation}
To simplify the notation, we will note in the sequel
$$\n {\cal E}*\big(P^{'}(\rho){\rm div}(\rho u)\big)=\n(\D)^{-1}\big(P^{'}(\rho){\rm div}(\rho u)\big).$$
\end{notation}
Finally we can now rewritte the system (\ref{0.6}) as follows:
\begin{equation}
\begin{cases}
\begin{aligned}
&\p_{t}a+(v_{1}+\frac{1}{\nu}v)\cdot\n a=(1+a){\rm div}(v_{1}+\frac{1}{\nu}v),\\
&\p_{t}v_{1}-(1+a){\cal A}v_{1}=f-u\cdot\n u+\frac{1}{\nu}\n(\D)^{-1}\big(P^{'}(\rho){\rm div}(\rho u)\big),\\
&a_{/ t=0}=a_{0},\;(v_{1})_{/ t=0}=(v_{1})_{0}.
\end{aligned}
\end{cases}
\label{0.7}
\end{equation}
where $v_{1}=u-\frac{1}{\nu}v$. In the sequel we will study this system by exctracting some uniform bounds in Besov spaces on
$(a,v_{1})$ as the in the following works \cite{AP}, \cite{H}. The advantage of the system (\ref{0.7}) is that we have \textit{kill} the coupling between $v_{1}$ and a term of pressure. Indeed in the works of R. Danchin \cite{DL}, \cite{DW}, the pressure was considered as a term of rest in the momentum equation, so it implied a strong relationship between the density and the velocity. In particular it was impossible to distinguish the index of integration for the Besov spaces.
\subsection{Proof of the existence}
\subsubsection*{Construction of approximate solutions}
We use a standard scheme:
\begin{enumerate}
\item We smooth out the data and get a sequence of smooth solutions $(a^{n},u^{n})_{n\in\mathbb{N}}$ to (\ref{0.6})
on a bounded interval $[0,T^{n}]$ which may depend on $n$. We set $v_{1}^{n}=u^{n}-v^{n}$ where ${\rm div}v^{n}=P(\rho^{n})-P(\bar{\rho})$.
\item We exhibit a positive lower bound $T$ for $T^{n}$, and prove uniform estimates on $(a^{n},u^{n})$ in the space
$$E_{T}=\widetilde{C}_{T}(B^{\NN}_{p,1})\times\big(\widetilde{C}_{T}(B^{\frac{N}{p_{1}}-1}_{p_{1},1}+B^{\NN+1}_{p,1})\cap\widetilde{L}^{1}_{T}(
B^{\frac{N}{p_{1}}+1}_{p_{1},1}+B^{\NN+2}_{p,1})\big).$$
More precisely to get this bounds we will need to study the behavior of $(a^{n},v_{1}^{n})$.
\item We use compactness to prove that the sequence $(a^{n},u^{n})$ converges, up to extraction, to a solution of (\ref{0.7}).
\end{enumerate}
Througout the proof, we denote $\underline{\nu}=\underline{b}\min(\mu,\lambda+2\mu)$ and $\bar{\nu}=\mu+|\mu+\lambda|$, and we assume (with no loss of generality) that $f$ belongs to $\widetilde{L}^{1}_{T}(B^{\frac{N}{p_{1}}}_{p_{1},1})$.
\subsubsection*{First step}
We smooth out the data as follows:
$$a_{0}^{n}=S_{n}a_{0},\;\;u_{0}^{n}=S_{n}u_{0}\;\;\;\mbox{and}\;\;\;f^{n}=S_{n}f.$$
Note that we have:
$$\forall l\in\mathbb{Z},\;\;\|\D_{l}a^{n}_{0}\|_{L^{p}}\leq\|\D_{l}a_{0}\|_{L^{p}}\;\;\;\mbox{and}\;\;\;\|a^{n}_{0}\|
_{B^{\frac{N}{p}}_{p,\infty}}\leq \|a_{0}\|_{B^{\frac{N}{p}}_{p,\infty}},$$
and similar properties for $u_{0}^{n}$ and $f^{n}$, a fact which will be used repeatedly during the next
steps. Now, according \cite{DW}, one can solve (\ref{0.6}) with the smooth data $(a_{0}^{n},u_{0}^{n},f^{n})$.
We get a solution $(a^{n},u^{n})$ on a non trivial time interval $[0,T_{n}]$ such that:
\begin{equation}
\begin{aligned}
&a^{n}\in\widetilde{C}([0,T_{n}),B^{N}_{2,1})\;\;\mbox{and}\;\;u^{n}\in\widetilde{C}([0,T_{n}),B^{\N-1}_{2,1})\cap
\widetilde{L}^{1}_{T_{n}}
(B^{\N+1}_{2,1}).
\end{aligned}
\label{a26}
\end{equation}
\subsubsection*{Uniform bounds}
Let $T_{n}$ be the lifespan of $(a_{n},u_{n})$, that is the supremum of all $T>0$ such that (\ref{0.1}) with initial data
$(a_{0}^{n},u_{0}^{n})$ has a solution which satisfies (\ref{a26}). Let $T$ be in $(0,T_{n})$.
We aim at getting uniform estimates in $E_{T}$ for $T$ small enough. For that, we need to introduce the solution $u^{n}_{L}$ to the linear system:
$$\p_{t}u_{L}^{n}-{\cal A}u_{L}^{n}=f^{n},\;\;u^{n}_{L}(0)=u^{n}_{0}-\frac{1}{\nu}\widetilde{v}_{0}\widetilde{v}_{0}.$$
Now, we set $\tilde{u}^{n}=u^{n}-u^{n}_{L}$ and the vectorfield $\widetilde{v}^{n}_{1}=\widetilde{u}^{n}-\frac{1}{\nu}\widetilde{v}^{n}$ with ${\rm div}\widetilde{v}^{n}=P(\rho^{n})$. We can check that $\widetilde{v}^{n}_{1}$ satisfies the parabolic system:
\begin{equation}
\begin{cases}
\begin{aligned}
&\p_{t}\widetilde{v}_{1}^{n}+(u_{L}^{n}+\frac{1}{\nu}\widetilde{v}^{n})\cdot\n \tilde{v}_{1}^{n}+\widetilde{v}_{1}^{n}\cdot\n u^{n}-(1+a^{n}){\cal A}\widetilde{v}_{1}^{n}=a^{n}
{\cal A}u_{L}^{n}-\frac{1}{\nu}(u_{L}^{n}\cdot\n \tilde{v}^{n}\\
&\hspace{3cm}+\frac{1}{\nu}\widetilde{v}^{n}\cdot\n \widetilde{v}^{n})-u_{L}^{n}\cdot\n u_{L}^{n}+\frac{1}{\nu}\n(\D)^{-1}(P^{'}(\rho^{n}){\rm div}(\rho^{n}u^{n})),\\
&(\widetilde{v}_{1}^{n})_{\ t=0}=0.
\end{aligned}
\end{cases}
\label{systemessen}
\end{equation}
which has been studied in proposition \ref{linearise}. 
Define $m\in\mathbb{Z}$ by:
\begin{equation}
m=\inf\{ p\in\mathbb{Z}/\;2\bar{\nu}\sum_{l\geq p}2^{l\frac{N}{p}}\|\D_{l}a_{0}\|_{L^{p}}\leq c\bar{\nu}\}
\label{def}
\end{equation}
where $c$ is small enough positive constant (depending only $N$) to be fixed hereafter. In the sequel we will need of a control on $a-S_{m}a$ small to apply proposition \ref{linearise}, so here $m$ is enough big (we explain how in the sequel).
Let:
$$\bar{b}=1+\sup_{x\in\R^{N}}a_{0}(x),\;A_{0}=1+2\|a_{0}\|_{B^{\NN}_{p,1}},\;U_{0}=\|u_{0}\|_{B^{\frac{N}{p_{1}}-1}_{p_{1},1}}+
\|a_{0}\|_{B^{\NN+1}_{p,1}}+
\|f\|_{L^{1}_{T}(B^{\frac{N}{p_{1}}-1}_{p_{1},1})},$$
and $\widetilde{U}_{0}=2CU_{0}+4C\bar{\nu}A_{0}$ (where $C^{'}$ is a constant embedding and
$C$ stands for a large enough constant depending only $N$ which will be determined when applying proposition \ref{produit1}, \ref{linearise} and \ref{transport1} in the following computations.) We assume that the following inequalities are fulfilled for some $\eta>0$:
$$
\begin{aligned}
&({\cal H}_{1})&\|a^{n}-S_{m}a^{n}\|_{\widetilde{L}^{\infty}_{T}(B^{\NN}_{p,1})}\leq c\underline{\nu}\bar{\nu}^{-1},\\
&({\cal H}_{2})&C\bar{\nu}^{2}T\|a^{n}\|^{2}_{\widetilde{L}^{\infty}_{T}(B^{\NN}_{p,1})}\leq 2^{-2m}\underline{\nu},\\
&({\cal H}_{3})&\frac{1}{2}\underline{b}\leq 1+a^{n}(t,x)\leq 2\bar{b}\;\;\mbox{for all}\;\;(t,x)\in[0,T]\times\R^{N},\\
&({\cal H}_{4})&\|a^{n}\|_{\widetilde{L}^{\infty}_{T}(B^{\NN}_{p,1})}\leq A_{0},
\end{aligned}
$$
$$
\begin{aligned}
&({\cal H}_{5})&\|u^{n}_{L}\|_{L^{1}_{T}(B^{\frac{N}{p_{1}}+1}_{p_{1},1}+B^{\NN+3}_{p,1})}\leq \eta,\\
&({\cal H}_{6})&\|\widetilde{v}_{1}^{n}\|_{\widetilde{L}^{\infty}_{T}(B^{\frac{N}{p_{1}}-1}_{p_{1},1}+B^{\NN+1}_{p,1})}+\underline{\nu}
\|\widetilde{v}_{1}^{n}\|_{L^{1}_{T}(B^{\frac{N}{p_{1}}+1}_{p_{1},1}+B^{\NN+2}_{p,1})}\leq \widetilde{U}_{0}\eta,\\
&({\cal H}_{7})&\|\widetilde{v}^{n}\|_{\widetilde{L}^{\infty}_{T}(B^{\frac{N}{p}+1}_{p,1})}\leq C^{'}A_{0},\\
&({\cal H}_{8})&\|\n u^{n}\|_{\widetilde{L}^{1}_{T}(B^{\frac{N}{p_{1}}}_{p_{1},1})+\widetilde{L}^{\infty}_{T}(B^{\NN}_{p,1})}\leq (\underline{\nu}^{-1}\widetilde{U}_{0}+1)\eta
\end{aligned}
$$
Remark that since:
$$1+S_{m}a^{n}=1+a^{n}+(S_{m}a^{n}-a^{n}),$$
assumptions $({\cal H}_{1})$ and $({\cal H}_{3})$ combined with the embedding $B^{\NN}_{p,1}\hookrightarrow L^{\infty}$ insure that:
\begin{equation}
\inf_{(t,x)\in[0,T]\times\R^{N}}(1+S_{m}a^{n})(t,x)\geq\frac{1}{4}\underline{b},
\label{inemin}
\end{equation}
provided $c$ has been chosen small enough (note that $\frac{\underline{\nu}}{\bar{\nu}}\leq\bar{b}$).\\
We are going to prove that under suitable assumptions on $T$ and $\eta$ (to be specified below) if condition $({\cal H}_{1})$ to $({\cal H}_{7})$ are satisfied, then they are actually satisfied with strict inequalities. Since all those conditions depend continuously on the time variable and
are strictly satisfied initially, a basic boobstrap argument insures that $({\cal H}_{1})$ to $({\cal H}_{8})$ are indeed satisfied for $T$.
First we shall assume that $\eta$ and $T$ satisfies:
\begin{equation}
C(1+\underline{\nu}^{-1}\widetilde{U}_{0})\eta+\frac{C^{'}}{\nu}A_{0}T
<\log 2
\label{1conduti}
\end{equation}
so that denoting $\widetilde{V}_{1}^{n}(t)=\int^{t}_{0}\|\n \widetilde{v}_{1}^{n}\|_{B^{\frac{N}{p_{1}}}_{p_{1},1}+B^{\NN+1}_{p,1}}d\tau$, $\widetilde{V}^{n}(t)=\frac{1}{\nu}\int^{t}_{0}\|\n \widetilde{v}^{n}\|_{B^{\frac{N}{p}}_{p,1}}d\tau$ and $U^{n}_{L}(t)=\int^{t}_{0}\|\n u^{n}_{L}\|_{B^{\frac{N}{p_{1}}+1}_{p_{1},1}+B^{\NN+3}_{p,1}}d\tau$, we have, according to $({\cal H}_{5})$ and $({\cal H}_{6})$:
\begin{equation}
e^{C(U^{n}_{L}+\widetilde{V}_{1}^{n}+\widetilde{V}^{n})(T)}<2\;\;\mbox{and}\;\;e^{C(U^{n}_{L}+\widetilde{V}_{1}^{n}+\widetilde{V}^{n})(T)}-1
\leq1.
\label{1ineimpca}
\end{equation}
In order to bound $a^{n}$ in $\widetilde{L}^{\infty}_{T}(B^{\NN}_{p,1})$, we apply inequality (\ref{22}) and get:
\begin{equation}
\|a^{n}\|_{\widetilde{L}^{\infty}_{T}(B^{\NN}_{p,1})}<1+2\|a_{0}\|_{B^{\NN}_{p,1}}=A_{0}.
\label{inetranspr}
\end{equation}
Hence $({\cal H}_{4})$ is satisfied with a strict inequality. $({\cal H}_{7})$ verifies a strict inequality, it follows from proposition \ref{singuliere} and $({\cal H}_{4})$.
Next, applying proposition \ref{chaleur} and proposition \ref{singuliere} yields:
\begin{equation}
\|u^{n}_{L}\|_{\widetilde{L}^{\infty}_{T}(B^{\frac{N}{p_{1}}-1}_{p_{1},1}+B^{\NN+1}_{p,1})}\leq U_{0},
\label{34}
\end{equation}
\begin{equation}
\begin{aligned}
&\kappa\nu\|u^{n}_{L}\|_{L^{1}_{T}(B^{\frac{N}{p_{1}}+1}_{p_{1},1}+B^{\NN+3}_{p,1})}\leq\sum_{l\in\mathbb{Z}}2^{l(\frac{N}{p_{1}}-1)}(1-e^{-\kappa\nu2^{2l}T})(\|\D_{l}u_{0}\|_{L^{p_{1}}}+\\
&\hspace{3cm}\|\D_{l}f\|_{L^{1}(\R^{+},L^{p_{1}})})+\leq\sum_{l\in\mathbb{Z}}2^{l(\frac{N}{p}+1)}(1-e^{-\kappa\nu2^{2l}T})
\|\D_{l}a_{0}\|_{L^{p}}.
\end{aligned}
\label{35}
\end{equation}
Hence taking $T$ such that:
\begin{equation}
\begin{aligned}
&\sum_{l\in\mathbb{Z}}2^{l(\frac{N}{p_{1}}-1)}(1-e^{-\kappa\nu2^{2l}T})(\|\D_{l}u_{0}\|_{L^{p_{1}}}
+\|\D_{l}f\|_{L^{1}(\R^{+},L^{p_{1}})})\\
&\hspace{2cm}+\leq\sum_{l\in\mathbb{Z}}2^{l(\frac{N}{p}+1)}(1-e^{-\kappa\nu2^{2l}T})
\|\D_{l}a_{0}\|_{L^{p}}<\kappa\eta\nu,
\end{aligned}
\label{36}
\end{equation}
insures that $({\cal H}_{5})$ is strictly verified.
Since $({\cal H}_{1})$, $({\cal H}_{2})$, $({\cal H}_{5})$, $({\cal H}_{6})$, $({\cal H}_{7})$ and (\ref{inemin}) are satisfied, proposition \ref{linearise} may be applied, we obtain:
$$
\begin{aligned}
&\|\widetilde{v}_{1}^{n}\|_{\widetilde{L}^{\infty}_{T}(B^{\frac{N}{p_{1}}-1}_{p_{1},1}+B^{\NN+1}_{p,1})}+\underline{\nu}
\|\widetilde{v}_{1}^{n}\|_{L^{1}_{T}(B^{\frac{N}{p_{1}}+1}_{p_{1},1}+B^{\NN+2}_{p,1})}\\
&\hspace{1cm}\leq Ce^{C(2U^{n}_{L}+2\widetilde{V}^{n}+\widetilde{V}_{1}^{n})(T)}\int^{T}_{0}\big(\|a^{n}{\cal A}u^{n}_{L}\|_{B^{\frac{N}{p_{1}}-1}_{p_{1},1}+B^{\NN}_{p,1}}
+\|u^{n}_{L}\cdot\n u^{n}_{L}\|_{B^{\frac{N}{p_{1}}-1}_{p_{1},1}+B^{\NN}_{p,1}}\\
&\hspace{1cm}+\|u_{L}^{n}\cdot\n \tilde{v}^{n}\|_{B^{\frac{N}{p_{1}}-1}_{p_{1},1}+B^{\NN}_{p,1}}+\|\widetilde{v}^{n}\cdot\n \widetilde{v}^{n}\|_{B^{\frac{N}{p}}_{p,1}}
+\|\n(\D)^{-1}(P^{'}(\rho^{n}){\rm div}(\rho^{n}u^{n}))\|_{B^{\frac{N}{p}}_{p,1}}\big) dt.
\end{aligned}
$$
As $\frac{N}{p}+\frac{N}{p_{1}}-1\geq 0$, $2\frac{N}{p}-1>0$ and by taking advantage of proposition \ref{produit1}, \ref{interpolation} and \ref{singuliere}, we get:
$$
\begin{aligned}
&\|\n(\D)^{-1}(P^{'}(\rho^{n}){\rm div}(\rho^{n}u^{n}))\|_{\widetilde{L}^{1}_{T}(B^{\frac{N}{p}}_{p,1})}
\leq C_{P}(1+\|a^{n}\|_{\widetilde{L}^{\infty}(B^{\frac{N}{p}}_{p,1})})(\sqrt{T}
\|\widetilde{v}_{1}^{n}\|_{\widetilde{L}^{2}_{T}(B^{\frac{N}{p_{1}}}_{p_{1},1}+B^{\NN+1}_{p,1})}\\
&\hspace{6,8cm}+\sqrt{T}\|u_{L}^{n}\|_{\widetilde{L}^{2}_{T}(B^{\frac{N}{p_{1}}}_{p_{1},1}+B^{\NN+1}_{p,1})}+T
\|a^{n}\|_{\widetilde{L}^{\infty}_{T}(B^{\frac{N}{p}}_{p,1})}),\\
&\|\widetilde{v}^{n}\cdot\n \widetilde{v}^{n}\|_{\widetilde{L}^{1}_{T}(B^{\frac{N}{p}}_{p,1})}\leq C_{1} T\|a^{n}\|^{2}_{\widetilde{L}_{T}^{\infty}(B^{\NN}_{p,1})}.
\end{aligned}
$$
We proceed similarly for the other terms and we end up with:
\begin{equation}
\begin{aligned}
&\|\widetilde{v}_{1}^{n}\|_{\widetilde{L}^{\infty}_{T}(B^{\frac{N}{p_{1}}-1}_{p_{1},1}+B^{\NN+1}_{p,1})}+\underline{\nu}
\|\widetilde{v}_{1}^{n}\|_{L^{1}_{T}(B^{\frac{N}{p_{1}}+1}_{p_{1},1}+B^{\NN+2}_{p,1})}\leq Ce^{C(2U^{n}_{L}+2\widetilde{V}^{n}+\widetilde{V}_{1}^{n})(T)}\\
&\times\biggl(C\|u^{n}_{L}\|_{L^{1}_{T}(B^{\frac{N}{p_{1}}+1}_{p_{1},1}+B^{\NN+3}_{p,1})}(\bar{\nu}\|a^{n}\|_{L^{\infty}_{T}
(B^{\frac{N}{p}}_{p,1})}
+\|u^{n}_{L}\|_{L^{\infty}_{T}(B^{\frac{N}{p_{1}}-1}_{p_{1},1}+B^{\NN+1}_{p,1})})+\\
&C_{1} T\|a^{n}\|^{2}_{\widetilde{L}_{T}^{\infty}(B^{\NN}_{p,1})}+C_{P}(1+\|a^{n}\|_{\widetilde{L}^{\infty}(B^{\frac{N}{p}}_{p,1})})(\sqrt{T}
\|\widetilde{v}_{1}^{n}\|_{\widetilde{L}^{2}_{T}(B^{\frac{N}{p_{1}}}_{p_{1},1}+B^{\NN+1}_{p,1})}\\
&+\sqrt{T}\|u_{L}^{n}\|_{\widetilde{L}^{2}_{T}(B^{\frac{N}{p_{1}}}_{p_{1},1}+B^{\NN+1}_{p,1})}+T
\|a^{n}\|_{\widetilde{L}^{\infty}_{T}(B^{\frac{N}{p}}_{p,1})})+T\|u^{n}_{L}\|_{L^{\infty}_{T}(B^{\frac{N}{p_{1}}-1}_{p_{1},1}+B^{\NN+1}_{p,1})}\times\\
&\hspace{10cm}\|a^{n}\|_{\widetilde{L}^{\infty}_{T}(B^{\frac{N}{p}}_{p,1})}\biggl).\\
\end{aligned}
\label{37}
\end{equation}
with $C=C(N)$, $C_{1}=C_{1}(N)$ and $C_{P}=(N,P,\underline{b},\bar{b})$. Now, using assumptions $({\cal H}_{4})$, $({\cal H}_{5})$ and $({\cal H}_{6})$, and inserting (\ref{1ineimpca}) in (\ref{37}) gives:
$$\|\widetilde{v}_{1}^{n}\|_{\widetilde{L}^{\infty}_{T}(B^{\frac{N}{p_{1}}-1}_{p_{1},1})}+
\|\widetilde{v}_{1}^{n}\|_{L^{1}_{T}(B^{\frac{N}{p_{1}}+1}_{p_{1},1})}\leq2C(\bar{\nu}A_{0}+U_{0})\eta+C_{1}TA_{0}(1+A_{0})+\sqrt{T}A_{0}U_{0},$$
hence $({\cal H}_{6})$ is satisfied with a strict inequality provided when $T$ verifies:
\begin{equation}
2C(\bar{\nu}A_{0}+U_{0})\eta+C_{1}TA_{0}(1+A_{0})+\sqrt{T}A_{0}U_{0}<C\bar{\nu}\eta.
\label{38}
\end{equation}
$({\cal H}_{8})$ verifies a strict inequality, it follows from proposition  $({\cal H}_{5})$, $({\cal H}_{6})$ and $({\cal H}_{7})$.
We now have to check whether $({\cal H}_{1})$ is satisfied with strict inequality. For that we apply proposition (\ref{transport2}) which yields for all $m\in\mathbb{Z}$,
\begin{equation}
\sum_{l\geq m}2^{l\N}\|\D_{l}a^{n}\|_{L^{\infty}_{T}(L^{p})}\leq\sum_{l\geq m}2^{l\NN}\|\D_{l}a_{0}\|_{L^{p}}+
(1+\|a_{0}\|_{B^{\N}_{p,1}})\big( e^{C(U^{n}_{L}+\widetilde{U}^{n})(T)}-1\big).
\label{39}
\end{equation}
Using (\ref{1conduti}) and $({\cal H}_{5})$, $({\cal H}_{6})$, we thus get:
$$\|a^{n}-S_{m}a^{n}\|_{L^{\infty}_{T}(B^{\NN}_{p,1})}\leq\sum_{l\geq m}2^{l\NN}\|\D_{l}a_{0}\|_{L^{p}}+\frac{C}{\log2}
(1+\|a_{0}\|_{B^{\NN}_{p,1}})(1+\underline{\nu}^{-1}\widetilde{L}_{0})\eta.$$
Hence $({\cal H}_{1})$ is strictly satisfied provided that $\eta$ further satisfies:
\begin{equation}
\frac{C}{\log2}(1+\|a_{0}\|_{B^{\NN}_{p,1}})(1+\underline{\nu}^{-1}\widetilde{U}_{0})\eta<\frac{c\underline{\nu}}{2\bar{\nu}}.
\label{40}
\end{equation}
In order to check whether $({\cal H}_{3})$ is satisfied, we use the fact that:
$$a^{n}-a_{0}=S_{m}(a^{n}-a_{0})+(Id-S_{m})(a^{n}-a_{0})+\sum_{l>n}\D_{l}a_{0},$$
whence, using $B^{\NN}_{p,1}\hookrightarrow L^{\infty}$ and assuming (with no loss of generality) that $n\geq m$,
$$
\begin{aligned}
&\|a^{n}-a_{0}\|_{L^{\infty}((0,T)\times\R^{N})}\leq C\big(\|S_{m}(a^{n}-a_{0})\|_{L^{\infty}_{T}(B^{\NN}_{p,1})}+
\|a^{n}-S_{m}a^{n}\|_{L^{\infty}_{T}(B^{\NN}_{p,1})}\\
&\hspace{9cm}+2\sum_{l\geq m}2^{l\NN}\|\D_{l}a_{0}\|_{L^{p}}\big).
\end{aligned}
$$
Changing the constant $c$ in the definition of $m$ and in (\ref{40}) if necessary, one can, in view of the previous
computations, assume that:
$$C\big(\|a^{n}-S_{m}a^{n}\|_{L^{\infty}_{T}(B^{\NN}_{p,1})}+2\sum_{l\geq m}2^{l\NN}\|\D_{l}a_{0}\|_{L^{p}}\big)\leq\frac{\underline{b}}{4}.$$
As for the term $\|S_{m}(a^{n}-a_{0})\|_{L^{\infty}_{T}(B^{\NN}_{p,1})}$, it may be bounded according proposition \ref{transport2}:
$$
\begin{aligned}
&\|S_{m}(a^{n}-a_{0})\|_{L^{\infty}_{T}(B^{\NN}_{p,1})}\leq(1+\|a_{0}\|_{B^{\NN}_{p,1}})(e^{C(\widetilde{V}_{1}^{n}+\widetilde{V}^{n}+U^{n}_{L})(T)}
-1)+C2^{2m}\sqrt{T}\|a_{0}\|_{B^{\NN}_{p,1}}\\
&\hspace{10cm}\times\|u^{n}\|_{L^{2}_{T}(B^{\frac{N}{p_{1}}}_{p_{1},1}+B^{\NN}_{p,1})}.
\end{aligned}
$$
Note that under assumptions $({\cal H}_{5})$, $({\cal H}_{6})$, (\ref{1conduti}) and (\ref{40}) ( and changing $c$ if necessary), the first term in the right-hand side may be bounded by $\frac{\underline{b}}{8}$.
Hence using interpolation, (\ref{34}) and the assumptions (\ref{1conduti}) and (\ref{40}), we end up with:
$$\|S_{m}(a^{n}-a_{0})\|_{L^{\infty}_{T}(B^{\NN}_{p,1})}\leq\frac{\underline{b}}{8}+C2^{m}\sqrt{T}\|a_{0}\|_{B^{\N}_{2,1}}
\sqrt{\eta(U_{0}+\widetilde{U}_{0}\eta)(1+\underline{\nu}^{-1}\widetilde{U}_{0}}.$$
Assuming in addition that $T$ satisfies:
\begin{equation}
C2^{m}\sqrt{T}\|a_{0}\|_{B^{\NN}_{p,1}}
\sqrt{\eta(U_{0}+\widetilde{U}_{0}\eta)(1+\underline{\nu}^{-1}\widetilde{U}_{0}}<\frac{\underline{b}}{8},
\label{42}
\end{equation}
and using the assumption $\underline{b}\leq1+a_{0}\leq\bar{b}$ yields $({\cal H}_{3})$ with a strict inequality.\\
One can now conclude that if $T<T^{n}$ has been chosen so that conditions (\ref{36}), (\ref{38}) and (\ref{42}) are satisfied (with $\eta$ verifying (\ref{1conduti}) and (\ref{40}), and $m$ defined in (\ref{def})
and $n\geq m$ then $(a^{n},u^{n})$ satisfies $({\cal H}_{1})$ to $({\cal H}_{8})$, thus is bounded independently of $n$
on $[0,T]$.\\
We still have to state that $T^{n}$ may be bounded by below by the supremum $\bar{T}$ of all times $T$
such that (\ref{36}), (\ref{38}) and (\ref{42}) are satisfied. This is actually a consequence of the uniform bounds we have just obtained, and of remark \ref{remark6} and proposition \ref{transport1}. Indeed, by combining all these informations, one can prove that if $T^{n}<\bar{T}$ then $(a^{n},u^{n})$ is actually in:
$$\widetilde{L}^{\infty}_{T^{n}}(B^{\N}_{2,1}\cap B^{\NN}_{p,1})\times\biggl(\widetilde{L}^{\infty}_{T^{n}}\big(B^{\N}_{2,1}\cap (B^{\frac{N}{p_{1}}-1}_{p_{1},1}+B^{\NN+1}_{p,1})\big)\cap L^{1}_{T^{n}}(B^{\N+1}_{2,1}\cap (B^{\frac{N}{p_{1}}-1}_{p_{1},1}+B^{\NN+2}_{p,1})\biggl)^{N}$$
hence may be continued beyond $\bar{T}$ (see the remark on the lifespan following the statement
in \cite{DL}). We thus have $T^{n}\geq\bar{T}$.
\subsubsection*{Compactness arguments}
We now have to prove that $(a^{n},u^{n})_{n\in\mathbb{N}}$ tends (up to a subsequence) to some function $(a,u)$ which belongs to $E_{T}$. Here we recall that:
$$E_{T}=\widetilde{C}([0,T],B^{\NN}_{p,1})\times\big(\widetilde{L}^{\infty}(B^{\frac{N}{p_{1}}-1}_{p_{1},1}+B^{\NN+1}_{p,1})\cap \widetilde{L}^{1}(
B^{\frac{N}{p_{1}}+1}_{p_{1},1}+B^{\NN+2}_{p,1})\big).$$
The proof is based on Ascoli's theorem and compact embedding for Besov spaces. As similar arguments have been employed in \cite{DL} or \cite{DW}, we only give the outlines of the proof.
\begin{itemize}
\item Convergence of $(a^{n})_{n\in\mathbb{N}}$:\\
We use the fact that $\widetilde{a}^{n}=a^{n}-a^{n}_{0}$ satisfies:
$$\p_{t}\widetilde{a}^{n}=-u^{n}\cdot\n a^{n}+(1+a^{n}){\rm div}u^{n}.$$
Since $(u^{n})_{n\in\mathbb{N}}$ is uniformly bounded in $\widetilde{L}^{1}_{T}(B^{\frac{N}{p_{1}}+1}_{p_{1},1}+B^{\NN+1}_{p,1})\cap L^{\infty}_{T}(B^{\frac{N}{p_{1}}-1}_{p_{1},1}+B^{\NN+1}_{p,1})$, it is by interpolation and the fact that $p_{1}\leq p$, also bounded in $L^{r}_{T}(B^{\frac{N}{p}-1+\frac{2}{r}}_{p,1})$ for any $r\in[1,+\infty]$.
By using the standard product laws in Besov spaces, we thus easily gather  that $(\p_{t}\widetilde{a}^{n})$ is uniformly bounded in $\widetilde{L}^{2}_{T}(B^{\NN-1}_{p,1})$. Hence $(\widetilde{a}^{n})_{n\in\mathbb{N}}$ is bounded in $\widetilde{L}^{\infty}_{T}(B^{\NN-1}_{p,1}\cap B^{\NN}_{p,1})$
and equicontinuous on $[0,T]$ with values in $B^{\NN-1}_{p,1}$. Since the embedding $B^{\NN-1}_{p,1}\cap B^{\NN}_{p,1}$ is (locally) compact, and $(a_{0}^{n})_{n\in\mathbb{N}}$ tends to $a_{0}$ in $B^{\NN}_{p,1}$, we conclude that $(a^{n})_{n\in\mathbb{N}}$ tends (up to extraction) to some distribution $a$. Given that $(a^{n})_{n\in\mathbb{N}}$ is bounded in $\widetilde{L}^{\infty}_{T}(B^{\NN}_{p,1})$, we actually have $a\in\widetilde{L}^{\infty}_{T}(B^{\NN}_{p,1})$.
\item Convergence of $(u^{n}_{L})_{n\in\mathbb{N}}$:\\
From the definition of $u^{n}_{L}$ and proposition \ref{chaleur}, it is clear that $(u^{n}_{L})_{n\in\mathbb{N}}$ tends to solution $u_{L}$ to:
$$\p_{t}u_{L}-{\cal A}u_{l}=f,\;\;u_{L}(0)=u_{0}-\frac{1}{\nu}.$$
in $\widetilde{L}^{\infty}_{T}(B^{\frac{N}{p_{1}}-1}_{p_{1},1}+B^{\NN+1}_{p,1})\cap \widetilde{L}^{1}_{T}(B^{\frac{N}{p_{1}}+1}_{p_{1},1}+B^{\NN+3}_{p,1})$.
\item  Convergence of $(\widetilde{v}_{1}^{n})_{n\in\mathbb{N}}$:\\
We use the fact that:
$$
\begin{aligned}
&\p_{t}\widetilde{v}_{1}^{n}=-(u_{L}^{n}+\frac{1}{\nu}\widetilde{v}^{n})\cdot\n \tilde{v}_{1}^{n}-\widetilde{v}_{1}^{n}\cdot\n u^{n}-\frac{1}{\nu}(u_{L}^{n}\cdot\n \tilde{v}^{n}-\frac{1}{\nu}\widetilde{v}^{n}\cdot\n \widetilde{v}^{n})+(1+a^{n}){\cal A}\widetilde{v}_{1}^{n}\\
&\hspace{4,5cm}+a^{n}
{\cal A}u_{L}^{n}-u_{L}^{n}\cdot\n u_{L}^{n}+\frac{1}{\nu}\n(\D)^{-1}(P^{'}(\rho^{n}){\rm div}(\rho^{n}u^{n})),\\
\end{aligned}
$$
As $(a^{n})_{n\in\mathbb{N}}$ is uniformly bounded in $L^{\infty}_{T}(B^{\NN}_{p,1})$ and $(u^{n})_{n\in\mathbb{N}}$ is uniformly bounded in $L^{\infty}_{T}(B^{\frac{N}{p_{1}}-1}_{p_{1},1}+B^{\NN+1}_{p,1})\cap L^{1}(B^{\frac{N}{p_{1}}+1}_{p_{1},1}+B^{\NN+1}_{p,1})$, it is easy to see that the the right-hand side is uniformly bounded in $\widetilde{L}^{\frac{4}{3}}_{T}(B^{\frac{N}{p_{1}}--\frac{3}{2}}_{p_{1},1})+\widetilde{L}^{\infty}(B^{\NN-1}_{p,1})$.
Hence $(\widetilde{v}_{1}^{n})_{n\in\mathbb{N}}$ is bounded in $\widetilde{L}^{\infty}_{T}(B^{\frac{N}{p_{1}}-1}_{p_{1},1}+B^{\NN+1}_{p,1})$ and equicontinuous on $[0,T]$ with values in $B^{\frac{N}{p_{1}}-1}_{p_{1},1}+B^{\frac{N}{p_{1}}-\frac{3}{2}}_{p_{1},1}$. This enables to conclude that
$(\widetilde{v}_{1}^{n})_{n\in\mathbb{N}}$ converges (up to extraction) to some function $\widetilde{v}_{1}\in
\widetilde{L}^{\infty}_{T}(B^{\frac{N}{p_{1}}-1}_{p_{1},1}+B^{\NN+1}_{p,1})\cap L^{1}_{T}(B^{\frac{N}{p_{1}}+1}_{p_{1},1}+B^{\NN+2}_{p,1})$.
\end{itemize}
By interpolating with the bounds provided by the previous step, one obtains better results of convergence so that one can pass to the limit in the mass equation and in the momentum equation. Finally by setting $u=\widetilde{v}_{1}+\widetilde{v}+u_{L}$, we conclude that
$(a,u)$ satisfies (\ref{0.6}).\\
In order to prove continuity in time for $a$ it suffices to make use of proposition \ref{transport1}. Indeed, $a_{0}$ is in $B^{\NN}_{p,1}$, and having $a\in \widetilde{L}^{\infty}_{T}(B^{\NN}_{p,1})$ and $u\in \widetilde{L}^{1}_{T}(B^{\frac{N}{p_{1}}+1}_{p_{1},1}+B^{\NN+1}_{p,1})$ insure that $\p_{t}a+u\cdot\n a$ belongs to $\widetilde{L}^{1}_{T}(B^{\NN}_{p,1})$. Similarly, continuity for $u$ may be proved by using that $(\widetilde{v}_{1})_{0}\in B^{\frac{N}{p_{1}}-1}_{p_{1},1}$ and that $(\p_{t}v_{1}-\mu\D v_{1})\in \widetilde{L}^{1}_{T}(B^{\frac{N}{p_{1}}-1}_{p_{1},1}+B^{\NN}_{p,1})$.
We conclude by using the fact that $u=v_{1}+\frac{1}{\nu}v$.
\subsection{The proof of the uniqueness}
\subsubsection*{Uniqueness when $1\leq p_{1}<2N$, $\frac{2}{N}<\frac{1}{p}+\frac{1}{p_{1}}$ and $N\geq 3$}
In this section, we focus on the cases $1\leq p_{1}<2N$, $\frac{2}{N}<\frac{1}{p}+\frac{1}{p_{1}}$, $N\geq 3$ and postpone the analysis of the other cases
(which turns out to be critical) to the next section.
Throughout the proof, we assume that we are given two solutions $(a^{1},u^{1})$ and $(a^{2},u^{2})$ of (\ref{0.6}). In the sequel we will show that $a^{1}=a^{2}$ and $v_{1}^{1}=v_{1}^{2}$ where $u^{i}=v_{1}^{i}+\widetilde{v}^{i}$. It will imply that $u^{1}=u^{2}$). We know that  $(a^{1},v_{1}^{1})$ and $(a^{2},v_{1}^{2})$ belongs to:
$$ \widetilde{C}([0,T]; B^{\NN}_{p,1})\times\big(\widetilde{C}([0,T];B^{\frac{N}{p_{1}}-1}_{p_{1},1}+B^{\NN+1}_{p,1})\cap \widetilde{L}^{1}(0,T;B^{\frac{N}{p_{1}}+1}_{p_{1},1}+B^{\NN+2}_{p,1})\big)^{N}.$$
Let $\delta a=a^{2}-a^{1}$, $\delta v=\widetilde{v}^{2}-\widetilde{v}^{1}$
and $\delta v_{1}=v_{1}^{2}-v_{1}^{1}$. The system for $(\delta a,\delta v_{1})$ reads:
\begin{equation}
\begin{cases}
\begin{aligned}
&\p_{t}\delta a+u^{2}\cdot\n\delta a=\delta a{\rm div} u^{2}+(\delta v_{1}+\frac{1}{\nu}\delta v)\cdot\n a^{1}+(1+a^{1}){\rm div}(\delta v_{1}+\frac{1}{\nu}\delta v),\\
&\p_{t}\delta v_{1}+u^{2}\cdot\delta \n v_{1}+\delta v_{1}\cdot\n u^{1}-(1+a^{1}){\cal A}\delta v_{1}=\delta a{\cal A}v_{1}^{2}-\frac{1}{\nu}(u^{2}\cdot\n\delta\widetilde{v}\\
&-\delta \widetilde{v}\cdot\n u^{1})+\n (\D)^{-1}\biggl((P^{'}(\rho^{2})-P^{'}(\rho^{1})){\rm div}(\rho^{2}u^{2})+P^{'}(\rho^{1}){\rm div}(\rho^{1}\delta u)\\
&\hspace{7,5cm}+
P^{'}(\rho^{1}){\rm div}((\rho^{2}-\rho^{1})u^{2})\biggl).
\end{aligned}
\end{cases}
\label{systemeuni}
\end{equation}
The function $\delta a$ may be estimated by taking advantage of proposition \ref{transport1} with $s=\NN-1$.
Denoting $U^{i}(t)=\|\n u^{i}\|_{\widetilde{L}^{1}(B^{\frac{N}{p_{1}}+1}_{p_{1},1}+B^{\NN+1}_{p,1})}$ for $i=1,2$, we get for all $t\in[0,T]$,
$$
\begin{aligned}
&\|\delta a(t)\|_{B^{\NN-1}_{p,1}}\leq C e^{C U^{2}(t)}\int^{t}_{0}e^{-CU^{2}(\tau)}\|\delta a{\rm div} u^{2}+(\delta v_{1}+\frac{1}{\nu}\delta v)\cdot\n a^{1}\\
&\hspace{7cm}+(1+a^{1}){\rm div}(\delta v_{1}+\frac{1}{\nu}\delta v)\|_{B^{\NN-1}_{p,1}}d\tau,
\end{aligned}
$$
Next using proposition \ref{produit1} and \ref{singuliere} we obtain:
$$
\begin{aligned}
&\|\delta a(t)\|_{B^{\NN-1}_{p,1}}\leq C e^{C U^{2}(t)}\int^{t}_{0}e^{-CU^{2}(\tau)}\|\delta a\|_{B^{\NN-1}_{p,1}}\big(\|u^{2}\|_{B^{\frac{N}{p_{1}}+1}_{p_{1},1}+B^{\NN+1}_{p,1}}+(1+2\|a_{1}\|_{B^{\NN}_{p,1}})\big)\\
&\hspace{7,7cm}+(1+2\|a_{1}\|_{B^{\NN}_{p,1}})\|\delta v_{1}\|_{B^{\frac{N}{p_{1}}}_{p_{1},1}+B^{\NN+1}_{p,1}}d\tau,
\end{aligned}
$$
Hence 
applying Gr\"onwall lemma, we get:
\begin{equation}
\|\delta a(t)\|_{B^{\NN-1}_{p,1}}\leq C e^{C U^{2}(t)}\int^{t}_{0}e^{-CU^{2}(\tau)}(1+\|a^{1}\|_{B^{\NN}_{p,1}})
\|\delta v_{1}\|_{B^{\frac{N}{p_{1}}}_{p_{1},1}+B^{\NN+1}_{p,1}}d\tau.
\label{ineunia}
\end{equation}
For bounding $\delta v_{1}$, we aim at applying proposition \ref{linearise} to the second equation of (\ref{systemeuni}).
So let us fix an integer $m$ such that:
\begin{equation}
1+\inf_{(t,x)\in[0,T]\times\R^{N}}S_{m}a^{1}\geq\frac{\underline{b}}{2}\;\;\mbox{and}\;\;\|a^{1}-
S_{m}a^{1}\|_{L^{\infty}_{T}(B^{\NN}_{p,1})}\leq c\frac{\underline{\nu}}{\bar{\nu}}.
\label{ineunicondi}
\end{equation}
Note since $a^{1}$ satisfies a transport equation with right-hand side in $\widetilde{L}^{1}_{T}(B^{\NN-1}_{p,1})$,
proposition \ref{transport1} guarantees that $a^{1}$ is in $\widetilde{C}_{T}(B^{\NN}_{p,1})$.
Hence such an integer does exist (see remark \ref{remark5}).
Now applying corollary \ref{linearise1} with $s=\frac{N}{p_{1}}-2$ and $s^{'}=\NN-1$ insures that for all time $t\in[0,T]$, we have:
$$
\begin{aligned}
&\|\delta v_{1}\|_{L^{1}_{t}(B^{\frac{N}{p_{1}}}_{p_{1},1}+B^{\NN+1}_{p,1})}\leq C e^{C U(t)}\int^{t}_{0}e^{-CU(\tau)}\big(\|\delta a{\cal A}v_{1}^{2}
-\frac{1}{\nu}(\delta v\cdot\n v_{1}^{1}+v_{1}^{1}\cdot\n\delta v)\\
&\hspace{7cm}-\frac{1}{\nu^{2}}(v^{1}\cdot\n\delta v+\delta v\cdot\n v^{2})\|_{B^{\frac{N}{p_{1}}-2}_{p_{1},1}+B^{\NN-1}_{p,1}}\big)d\tau,
\end{aligned}
$$
with $U(t)=U^{1}(t)+U^{2}(t)+2^{2m}\underline{\nu}^{-1}\bar{\nu}^{2}\int^{t}_{0}\|a^{1}\|^{2}_{B^{\NN}_{p,1}}d\tau$.\\
Hence, applying proposition \ref{produit1} we get:
\begin{equation}
\begin{aligned}
&\|\delta v_{1}\|_{\widetilde{L}^{1}_{t}(B^{\frac{N}{p_{1}}}_{p_{1},1}+B^{\NN+1}_{p,1})}\leq C e^{C U(t)}\int^{t}_{0}e^{-CU(\tau)}\big(1+\|a^{1}\|_{B^{\NN}_{p,1}}
+\|a^{2}\|_{B^{\NN}_{p,1}}\\
&\hspace{7cm}+\|v_{1}^{2}\|_{B^{\frac{N}{p_{1}}+1}_{p_{1},1}+B^{\NN+2}_{p,1}}
\big)\|\delta a\|_{B^{\NN-1}_{p,1}}d\tau.
\end{aligned}
\label{ineunide}
\end{equation}
Finally plugging (\ref{ineunia}) in (\ref{ineunide}), we get for all $t\in[0,T_{1}]$,
$$
\begin{aligned}
&\|\delta v_{1}\|_{\widetilde{L}^{1}_{t}(B^{\frac{N}{p_{1}}}_{p_{1},1}+B^{\NN+1}_{p,1})}\leq C e^{C U(t)}\int^{t}_{0}\big(1+\|a^{1}\|_{B^{\NN}_{p,1}}
+\|a^{2}\|_{B^{\NN}_{p,1}}+\|v_{1}^{2}\|_{B^{\frac{N}{p_{1}}+1}_{p_{1},1}+B^{\NN+2}_{p,1}}
\big)\\
&\hspace{9,5cm}\times\|\delta v_{1}\|_{B^{\frac{N}{p_{1}}}_{p_{1},1}+B^{\NN+1}_{p,1}}d\tau.
\end{aligned}
$$
Since $a^{1}$ and $a^{2}$ are in $L^{\infty}(B^{\NN}_{p,1})$ and $v_{1}^{2}$ belongs to $L^{1}_{T}(B^{\frac{N}{p_{1}}+1}_{p_{1},1}+B^{\NN+2}_{p,1})$,
applying Gr\"onwall lemma yields $\delta v_{1}=0$, an $[0,T]$.
\subsubsection*{Uniqueness when:$\frac{2}{N}=\frac{1}{p_{1}}+\frac{1}{p}$ or $p_{1}=2N$ or $N=2$.}
The above proof fails in dimension two. One of the reasons why is that the product of functions does not map
$B^{\NN}_{p,1}\times B^{\frac{N}{p_{1}}-2}_{p_{1},1}$ in $B^{\frac{N}{p_{1}}-2}_{p_{1},1}$ but only in the larger space $B^{\frac{N}{p_{1}}-2}_{p_{1},\infty}$. This induces us to bound $\delta a$ in $\L_{T}^{\infty}(B^{\NN-1}_{p,\infty})$ and $\delta v_{1}$ in $L_{T}^{\infty}(B^{\frac{N}{p_{1}}-2}_{p_{1},\infty}+B^{\NN}_{p,\infty})\cap L^{1}_{T}(B^{\frac{N}{p_{1}}}_{p_{1},\infty}+B^{\NN+1}_{p,\infty})$ (or rather, in the widetilde version of those spaces, see below). Yet, we are in trouble because due to $B^{\frac{N}{p_{1}}}_{p_{1},\infty}$ is not embedded in $L^{\infty}$, the term $\delta v_{1}\cdot\n a^{1}$ in the right hand-side of the first equation of (\ref{systemeuni}) cannot be estimated properly. As noticed in \cite{DU}, this second difficulty may be overcome by making use of logarithmic interpolation and Osgood lemma ( a substitute for Gronwall inequality).
Let us now tackle the proof. Fix an integer $m$ such that:
\begin{equation}
1+\inf_{(t,x)\in[0,T]\times\R^{N}}S_{m}a^{1}\geq\frac{\underline{b}}{2}\;\;\mbox{and}\;\;\|a^{1}-
S_{m}a^{1}\|_{\widetilde{L}^{\infty}_{T}(B^{\NN}_{p,1})}\leq c\frac{\underline{\nu}}{\bar{\nu}},
\label{47}
\end{equation}
and define $T_{1}$ as the supremum of all positive times $t$ such that:
\begin{equation}
t\leq T\;\;\mbox{and}\;\;t\bar{\nu}^{2}\|a^{1}\|_{\widetilde{L}^{\infty}_{T}(B^{\NN}_{p,1})}\leq c2^{-2m}\underline{\nu}.
\label{48}
\end{equation}
Remark that the proposition \ref{transport1} ensures that $a^{1}$ belongs to $\widetilde{C}_{T}(B^{\NN}_{p,1})$
so that the above two assumptions are satisfied if $m$ has been chosen large enough.
For bounding $\delta a$ in $L^{\infty}_{T}(B^{\NN-1}_{p,\infty})$,
we apply proposition \ref{transport1} with $r=+\infty$ and $s=0$. We get (with the notation of the previous section):
$$
\begin{aligned}
&\forall t\in[0,T],\;\;\|\delta a(t)\|_{B^{\NN-1}_{p,\infty}}\leq Ce^{CU^{2}(t)}\int^{t}_{0}
e^{-CU^{2}(\tau)}\|\delta a{\rm div} u^{2}+(\delta v_{1}+\frac{1}{\nu}\delta v)\cdot\n a^{1}\\
&\hspace{7,5cm}+(1+a^{1}){\rm div}(\delta v_{1}+\frac{1}{\nu}\delta v)\|_{B^{\NN-1}_{p,\infty}}d\tau,
\end{aligned}
$$
hence using that the product of two functions maps $B^{\NN-1}_{p,\infty}\times
B^{\frac{N}{p_{1}}}_{p_{1},1}$ in $B^{\NN-1}_{p,\infty}$, and applying Gronwall lemma,
\begin{equation}
\|\delta a(t)\|_{B^{\NN-1}_{p,\infty}}\leq Ce^{CU^{2}(t)}\int^{t}_{0}
e^{-CU^{2}(\tau)}(1+\|a^{1}\|_{B^{\NN}_{p,1}})\|\delta v_{1}\|_{B^{\frac{N}{p_{1}}}_{p_{1},1}+B^{\NN+1}_{p,1}}d\tau.
\label{49}
\end{equation}
Next, using proposition \ref{linearise1} combined with proposition \ref{produit1} and corollary \ref{produit2} in order to bound the nonlinear terms, we get for all $t\in[0,T_{1}]$,:
\begin{equation}
\begin{aligned}
&\|\delta v_{1}\|_{\widetilde{L}^{1}_{T}(B^{\frac{N}{p_{1}}-2}_{p_{1},\infty}+B^{\NN+1}_{p,\infty})}\leq Ce^{C(U^{1}+U^{2})(t)}\int^{t}_{0}(1+\|a^{1}\|_{B^{\NN}_{p,1}}+\|a^{2}\|_{B^{\NN}_{p,1}}\\
&\hspace{6,5cm}+\|v_{1}^{2}\|_{B^{\frac{N}{p_{1}}+1}_{p_{1},1}+B^{\NN+2}_{p,1}})\|\delta a\|_{B^{\NN-1}_{p,\infty}}d\tau.
\end{aligned}
\label{50}
\end{equation}
In order to control the term $\|\delta v_{1}\|_{B^{\frac{N}{p_{1}}}_{p_{1},1}+B^{\NN+1}_{p,1}}$ which appears in the right-hand side of (\ref{49}), we make
use of the following logarithmic interpolation inequality whose proof may be found in \cite{DU}, page 120:
\begin{equation}
\begin{aligned}
&\|\delta v_{1}\|_{L^{1}_{t}(B^{\frac{N}{p_{1}}}_{p_{1},1}+B^{\NN+1}_{p,1})}\lesssim\\
&\|\delta v_{1}\|_{\widetilde{L}^{1}_{t}(B^{\frac{N}{p_{1}}}_{p_{1},\infty})}\log\big(
e+\frac{\|\delta v_{1}\|_{\widetilde{L}^{1}_{t}(B^{\frac{N}{p_{1}}-1}_{p_{1},\infty})}+\|\delta v_{1}\|_{\widetilde{L}^{1}_{t}(B^{\frac{N}{p_{1}}+1}_{p_{1},\infty})}}{\|\delta v_{1}\|_{\widetilde{L}^{1}_{t}(B^{\frac{N}{p_{1}}}_{p_{1},\infty})}}\big)\\
&\hspace{2cm}+\|\delta v_{1}\|_{\widetilde{L}^{1}_{t}(B^{\NN+1}_{p,\infty})}\log\big(
e+\frac{\|\delta v_{1}\|_{\widetilde{L}^{1}_{t}(B^{\NN}_{p,\infty})}+\|\delta v_{1}\|_{\widetilde{L}^{1}_{t}(B^{\NN+2}_{p,\infty})}}{\|\delta v_{1}\|_{\widetilde{L}^{1}_{t}(B^{\NN}_{p,\infty})}}\big).
\end{aligned}
\label{51}
\end{equation}
Because $v_{1}^{1}$ and $v_{2}^{2}$ belong to $\widetilde{L}^{\infty}_{T}(B^{\frac{N}{p_{1}}-1}_{p_{1},1}+B^{\NN+1}_{p,1})\cap L^{1}_{T}(B^{\frac{N}{p_{1}}+1}_{p_{1},1}+B^{\NN+2}_{p,1})$,
the numerator in the right-hand side may be bounded by some constant $C_{T}$
depending only on $T$ and on the norms of $v_{1}^{1}$ and $v_{1}^{2}$.
Therefore inserting (\ref{49}) in (\ref{50}) and taking advantage of (\ref{51}), we end up for
all $t\in[0,T_{1}]$ with:
$$
\begin{aligned}
&\|\delta v_{1}\|_{\widetilde{L}^{1}_{T}(B^{\frac{N}{p_{1}}}_{p_{1},1}+B^{\NN+1}_{p,1})}\leq C(1+\|a^{1}\|_{\widetilde{L}^
{\infty}_{T}(B^{\NN}_{p,1})})\\
&\hspace{0,5cm}\times\int^{t}_{0}(1+\|a^{1}\|_
{B^{\NN}_{p,1}}+\|a^{2}\|_{B^{\NN}_{p,1}}+\|v_{1}^{2}\|_{B^{\frac{N}{p_{1}}+1}_{p_{1},1}+B^{\NN+2}_{p,1}})\|\delta v_{1}\|_
{\widetilde{L}^{1}_{t}(B^{\frac{N}{p_{1}}}_{p_{1},\infty})}\\
&\hspace{6cm}\times\log(e+C_{T}\|\delta v_{1}\|^{-1}_
{\widetilde{L}^{1}_{\tau}(B^{\frac{N}{p_{1}}}_{p_{1},\infty}+B^{\NN+1}_{p,\infty})}\big)d\tau.
\end{aligned}
$$
Since the function $t\rightarrow\|a^{1}(t)\|_{B^{\NN}_{p,1}}+\|a^{2}(t)\|_{B^{\NN}_{p,1}}+\|v_{1}^{2}(t)\|_{B^{\frac{N}{p_{1}}+1}_{p_{1},1}+
B^{\NN+2}_{p,1}}$
is integrable on $[0,T]$, and:
$$\int^{1}_{0}\frac{dr}{r\log(e+C_{T}r^{-1})}=+\infty$$
Osgood lemma yields $\|\delta v_{1}\|_{\widetilde{L}^{1}_{T}(B^{\frac{N}{p_{1}}}_{p_{1},1}+B^{\NN+1}_{p,1})}=0$. Note that the
definition of $m$ depends only on $T$ and that (\ref{ineunicondi}) is satisfied on $[0,T]$.
Hence, the above arguments may be repeated on $[T_{1},2T_{1}]$, $[2T_{1},3T_{1}]$,etc.
until the whole interval $[0,T]$ is exhausted. This yields uniqueness on $[0,T]$ for $a$ and $v_{1}$
which implies uniqueness for $u$.
\subsection{Proof of corollary \ref{coro11}}
The proof follows the same line as theorem \ref{theo1} except concerning the term of rest $\n(\D)^{-1}(P^{'}(\rho){\rm div}(\rho u))$
in the momentum equation of system (\ref{0.7}). Indeed in our case this term can write simplify on the form $\rho u$. In this case we control this term in $\widetilde{L}^{2}(B^{\NN}{p,1})$ without imposing additional conditions on $p$ of type $2\NN-1>0$.\\
Now the difficulty is to prove the uniqueness. For that we use the main theorem of D. Hoff in \cite{5H5} which is a result weak-strong uniqueness. In this article, D. Hoff has two solutions $(\rho,u)$ and $(\rho_{1},u_{1})$ with the sme initial data $(\rho_{0},u_{0})$ and he show that under some hypothesis of regularity on $(\rho_{1},u_{1})$ and $(\rho_{2},u_{2})$ then $\rho_{1}=\rho_{2}$, $u_{1}=u_{2}$.
We now discuss that our solution check the conditions required in \cite{5H5}. More precisely we have to show that our solution $(\rho,u)$
verify all the hypothesis asked on $(\rho_{1},u_{1})$ and $(\rho_{2},u_{2})$. The check is easy and tedious, but only
one hypothesis required to be carreful and is in fact the main condition why D. Hoff does not get global strong solution in dimension $N=3$ for the solutions built in \cite{5H4}.
We need to check that $u\in L^{\infty}_{loc}((0,T],L^{\infty})$ and $\n u\in L^{1}((0,T),L^{\infty})$.
In our case we have $\n u=\n v_{1}+\frac{1}{\nu}\n v$ where we recall that ${\rm div}v=P(\rho)-P(\bar{\rho})$.
We know that by interpolation $\n v_{1}\in L^{1}_{T}(B^{\frac{N}{p_{1}}}_{p_{1},1}+B^{\NN+1}_{p,1})\h L^{1}_{T}(L^{\infty})$
and by proposition \ref{singuliere} $\n v\in L^{\infty}_{T}(B^{\NN}_{p,1})$. We obtain then $\n u\in L^{1}_{T}(L^{\infty})$.
We have now to show that $u\in L^{\infty}_{T}(L^{\infty})$. In fact we have just to apply classical energy inequalities, so we multiply the momentum equation by $u|u|^{p_{1}-2}$
$$
\begin{aligned}
&\frac{1}{p_{1}}\int_{\R^{N}}\rho|u|^{p_{1}}(t,x)dx+\mu\int^{t}_{0}|u|^{p_{1}-2}|\n u|^{2}(t,x)dtdx+\frac{p_{1}-2}{4}\mu\int^{t}_{0}
|u|^{p_{1}-4}|\n|u|^{2}|^{2}(t,x)dxdt\\
&\hspace{1,7cm}+\int^{t}_{0}\int_{\R^{N}} \big(P(\rho)-P(\bar{\rho})\big)\big({\rm div}u|u|^{p_{1}-2}+(p_{1}-2)
\sum_{i,k}u_{i}u_{k}\p_{i}u_{k}|u|^{p_{1}-4}\big)(t,x)dtdx\\
&\hspace{11cm}\leq\int_{\R^{N}}\rho_{0}|u_{0}|^{p_{1}}dx.
\end{aligned}
$$
By Young's inequalities and the fact that $P(\rho)-P(\bar{\rho})$ belongs in $L^{\infty}(L^{2}\cap L^{\infty})$ we obtain that for all $p_{1}\in [1,+\infty[$, $\rho^{\frac{1}{p_{1}}}u\in L^{\infty}(L^{p_{1}})$ and:
$$\|\rho^{\frac{1}{p_{1}}}u\|_{L^{\infty}(L^{p_{1}})}\leq C_{0},$$
where $C_{0}$ depend only of the initial data. As $\frac{1}{\rho}\in L^{\infty}$, we conclude that $u$ is uniformly bounded
in all spaces $L^{\infty}(L^{p_{1}})$ with $p_{1}\in [1,+\infty[$. We conclude then $u\in L^{\infty}_{T}(L^{\infty})$.
\section{Continuation criterions}
\label{section6}
\subsection*{Proof of theorem \ref{theo3}}
We now prove theorem \ref{theo3}. We have assumed here that $\rho_{0}^{\frac{1}{p_{1}}}u_{0}\in L^{p_{1}}$ with $p_{1}>N$. We want show that with our hypothesis in particular that $a\in L^{\infty}_{T}$ and $1+a$ bounded away on $[0,T]$, then we are able to show that
$\rho^{\frac{1}{p_{1}}}u\in L^{\infty}_{T}(L^{p_{1}})$. In this case as $1+a$ is bounded away we show that $u\in L^{\infty}_{T}(L^{p_{1}})$ and
 by embedding we get $u\in L^{\infty}_{T}(B^{\frac{N}{p_{1}}-1}_{p_{1},1})$ as $\frac{N}{p_{1}}\leq 0$. We can next conclude by the fact that $(a(T,\cdot),u(T,\cdot))\in (B^{\NN}_{p,1}\times B^{\frac{N}{p_{1}}-1}_{p_{1},1})$ so that we can extend our solutions.
Finally we have just to show that $\rho^{\frac{1}{p_{1}}}u\in L^{\infty}_{T}(L^{p_{1}})$, in this goal we have just to apply classical energy inequality. We multiply the momentum equation by $u|u|^{p_{1}-2}$ and we get after integration by part:
$$
\begin{aligned}
&\frac{1}{p_{1}}\int_{\R^{N}}\rho|u|^{p_{1}}(t,x)dx+\mu\int^{t}_{0}|u|^{p_{1}-2}|\n u|^{2}(t,x)dtdx+\frac{p_{1}-2}{4}\mu\int^{t}_{0}
|u|^{p_{1}-4}|\n|u|^{2}|^{2}(t,x)dxdt\\
&+\lambda\int^{t}_{0}\int_{\R^{N}}({\rm div}u)^{2}|u|^{p_{1}-2}(t,x)dtdx+\lambda\frac{p_{1}-2}{2}\int^{t}_{0}\int_{\R^{N}}
{\rm div}u\sum_{i}u_{i}\p_{i}|u|^{2}|u|^{p_{1}-4}(t,x)dtdx-\\
&\hspace{2cm}\int^{t}_{0}\int_{\R^{N}} \big(P(\rho)-P(\bar{\rho})\big)\big({\rm div}u|u|^{p_{1}-2}+(p_{1}-2)
\sum_{i,k}u_{i}u_{k}\p_{i}u_{k}|u|^{p_{1}-4}\big)(t,x)dtdx\\
&\hspace{11cm}\leq\int_{\R^{N}}\rho_{0}|u_{0}|^{p_{1}}dx.
\end{aligned}
$$
By Young's inequalities, inequality (\ref{inegaliteviscosite}) and the fact that $P(\rho)-P(\bar{\rho})$ belongs in $L^{\infty}(L^{1}\cap L^{\infty})$ we conclude the proof.
\section{Appendix}
This section is devoted to the proof of commutator estimates which have been used in section $2$ and $3$. They are based on
paradifferentiel calculus, a tool introduced by J.-M. Bony in \cite{BJM}. The basic idea of paradifferential calculus is that
any product of two distributions $u$ and $v$ can be formally decomposed into:
$$uv=T_{u}v+T_{v}u+R(u,v)=T_{u}v+T^{'}_{v}u$$
where the paraproduct operator is defined by $T_{u}v=\sum_{q}S_{q-1}u\D_{q}v$, the remainder operator $R$, by
$R(u,v)=\sum_{q}\D_{q}u(\D_{q-1}v+\D_{q}v+\D_{q+1}v)$ and $T^{'}_{v}u=T_{v}u+R(u,v)$.
Inequalities (\ref{12}) and (\ref{18}) are consequence of the following lemma:
\begin{lemme}
\label{alemme2}
Let $1\leq p_{1}\leq p\leq+\infty$ and $\sigma\in(-\min(\NN,\frac{N}{p_{1}^{'}}),\NN+1]$. There exists a sequence $c_{q}\in l^{1}(\mathbb{Z})$ such that $\|c_{q}\|_{l^{1}}=1$
and a constant
$C$ depending only on $N$ and $\sigma$ such that:
\begin{equation}
\forall q\in\mathbb{Z},\;\;\|[v\cdot\n,\D_{q}]a\|_{L^{p_{1}}}\leq C c_{q}2^{-q\sigma}\|\n v\|_{B^{\NN}_{p,1}}
\|a\|_{B^{\sigma}_{p_{1},1}}.
\label{52}
\end{equation}
In the limit case $\sigma=-\min(\NN,\frac{N}{p_{1}^{'}})$, we have:
\begin{equation}
\forall q\in\mathbb{Z},\;\;\|[v\cdot\n,\D_{q}]a\|_{L^{p_{1}}}\leq C c_{q}2^{q\NN}\|\n v\|_{B^{\NN}_{p,1}}
\|a\|_{B^{-\frac{N}{p_{1}}}_{p,\infty}}.
\label{53}
\end{equation}
Finally, for all $\sigma>0$ and $\frac{1}{p_{2}}=\frac{1}{p_{1}}-\frac{1}{p}$, there exists a constant $C$ depending only on $N$ and on $\sigma$ and a sequence
$c_{q}\in l^{1}(\mathbb{Z})$ with norm
$1$ such that:
\begin{equation}
\forall q\in\mathbb{Z},\;\;\|[v\cdot\n,\D_{q}]v\|_{L^{p}}\leq C c_{q}2^{-q\sigma}(\|\n v\|_{L^{\infty}}\|v\|_{B^{\sigma}_{p_{1},1}}+\|\n v\|_{L^{p_{2}}}\|\n v\|_{B^{\sigma-1}_{p,1}}).
\label{54}
\end{equation}
\end{lemme}
{\bf Proof:} These results are proved in \cite{BCD} chapter $2$.
{\hfill $\Box$}\\
Inequality (\ref{13}) is a consequence of the following lemma:
\begin{lemme}
\label{alemme3}
Let $1\leq p_{1}\leq p\leq+\infty$ and
$\alpha\in(1-\NN,1]$, $k\in\{1,\cdots,N\}$ and $R_{q}=\D_{q}(a\p_{k}w)-\p_{k}(a\D_{q}w)$. There
exists $c=c(\alpha,N,\sigma)$ such that:
\begin{equation}
\sum_{q}2^{q\sigma}\|R_{q}\|_{L^{p_{1}}}\leq C\|a\|_{B^{\NN+\alpha}_{p,1}}\|w\|_{B^{\sigma+1-\alpha}_{p_{1},1}}
\label{57}
\end{equation}
whenever $-\NN<\sigma\leq\alpha+\NN$.\\
In the limit case $\sigma=-\NN$, we have for some constant $C=C(\alpha,N)$:
\begin{equation}
\sup_{q}2^{-q\frac{N}{p}}\|R_{q}\|_{L^{p_{1}}}\leq C\|a\|_{B^{\NN+\alpha}_{p,1}}\|w\|_{B^{-\frac{N}{p_{1}}+1-\alpha}_{p_{1},\infty}}.
\label{58}
\end{equation}
\end{lemme}
{\bf Proof}
The proof is almost the same as the one of lemma A3 in \cite{DL}.
It is based on Bony's decomposition which enables us to split $R_{q}$
into:
$$R_{q}=\underbrace{\p_{k}[\D_{q},T_{a}]w}_{R_{q}^{1}}-\underbrace{\D_{q}T_{\p_{k}a}w}_{R_{q}^{2}}+\underbrace{\D_{q}T_{\p_{k}w}w}_{R_{q}^{3}}
+\underbrace{\D_{q}R(\p_{k}w,a)}_{R_{q}^{4}}-\underbrace{\p_{k}T^{'}_{\D_{q}w}a}_{R_{q}^{5}}.$$
Using the fact that:
$$R^{1}_{q}=\sum^{q+4}_{q^{'}=q-4}\p_{k}[\D_{q},S_{q^{'}-1}a]\D_{q^{'}}w,$$
and the mean value theorem, we readily get under the hypothesis that $\alpha\leq1$,
\begin{equation}
\sum_{q}2^{q\sigma}\|R^{1}_{q}\|_{L^{p_{1}}}\lesssim\|\n a\|_{B^{\alpha-1}_{\infty,1}}\|w\|_{B^{\sigma+1-\alpha}_{p_{1},1}}.
\label{59}
\end{equation}
Standard continuity results for the paraproduct insure that $R^{2}_{q}$ satisfies (\ref{59}) and that:
\begin{equation}
\sum_{q}2^{q\sigma}\|R^{1}_{q}\|_{L^{p_{1}}}\lesssim\|\n w\|_{B^{\sigma-\alpha-\frac{N}{p_{1}}}_{\infty,1}}\|a\|_{B^{\NN+\alpha}_{p,1}}.
\label{60}
\end{equation}
provided $\sigma-\alpha-\NN\leq0.$
Next, standard continuity result for the remainder insure that under the hypothesis $\sigma>-\NN$, we have:
\begin{equation}
\sum_{q}2^{q\sigma}\|R^{1}_{q}\|_{L^{p_{1}}}\lesssim\|\n w\|_{B^{\sigma-\alpha}_{p_{1},1}}\|a\|_{B^{\NN+\alpha}_{p,1}}.
\label{61}
\end{equation}
For bounding $R^{5}_{q}$ we use the decomposition:
$R^{5}_{q}=\sum_{q^{'}\geq q-3}\p_{k}(S_{q^{'}+2}\D_{q}w\D_{q^{'}}a),$
which leads (after a suitable use of Bernstein and H\"older inequalities) to:
$$2^{q\sigma}\|R^{5}_{q}\|_{L^{p_{1}}}\lesssim\sum_{q^{'}\geq q-2}2^{(q-q^{'})(\alpha+\frac{N}{p_{1}}-1)}2^{q(\sigma+1-\alpha)}
\|\D_{q}w\|_{L^{p_{1}}}2^{q^{'}(\NN+\alpha)}
\|\D_{q^{'}}a\|_{L^{p}}.$$
Hence, since $\alpha+\NN-1>0$, we have:
$$\sum_{q}2^{q\sigma}\|R^{5}_{q}\|_{L^{p}}\lesssim\|\n w\|_{B^{\sigma+1-\alpha}_{p_{1},1}}\|a\|_{B^{\NN+\alpha}_{p,1}}.$$
Combining this latter inequality with (\ref{59}), (\ref{60}) and (\ref{61}), and using the embedding $B^{\NN}_{p,1}\hookrightarrow
B^{r-\NN}_{\infty,1}$
for $r=\NN+\alpha-1$, $\sigma_\alpha$ completes the proof of (\ref{57}).\\
The proof of (\ref{58}) is almost the same: for bounding $R^{1}_{q}$, $R^{2}_{q}$, $R^{3}_{q}$ and $R^{5}_{q}$, it is just a matter of changing $\sum_{q}$ into
$\sup_{q}$. 
\null{\hfill $\Box$}
\begin{remarka}
For proving proposition \ref{linearise1}, we shall actually use the following non-stationary version of inequality (\ref{58}):
$$\sup_{q}2^{-q\NN}\|R_{q}\|_{L^{1}_{T}(L^{p_{1}})}\leq C\|a\|_{\widetilde{L}^{\infty}_{T}(B^{\frac{N}{p}+\alpha}_{p,1})}
\|w\|_{\widetilde{L}^{1}_{T}(B^{-\frac{N}{p_{1}}+1-\alpha}_{p_{1},\infty})},$$
which may be easily proved by following the computations of the previous proof, dealing with the time dependence according to H\"older inequality.
\label{remarque7}
\end{remarka}

\end{document}